\documentclass[11pt]{article}
\usepackage{amsmath}
\usepackage{amssymb}
\usepackage{graphicx}
\usepackage{xcolor}
\def\n{\noindent}
\def\L{{\bf L}}

\def\D{{\mathcal D}}
\def\P{{\mathcal P}}

\def\dint{\int\!\!\int}
\def\ve{\varepsilon}
\def\n{\noindent}
\def\vp{\varphi}

\def\Tilde{\widetilde}

\def\C{{\mathcal C}}

\def\N{{\mathcal N}}

\def\Q{{\mathcal Q}}

\def\ds{\displaystyle}
\def\sqr#1#2{\vbox{\hrule height .#2pt
\hbox{\vrule width .#2pt height #1pt \kern #1pt
\vrule width .#2pt}\hrule height .#2pt }}
\def\square{\sqr74}
\def\endproof{\hphantom{MM}\hfill\llap{$\square$}\goodbreak}

\def\bega{\begin{array}}
\def\enda{\end{array}}
\def\begi{\begin{itemize}}
\def\endi{\end{itemize}}

\def\K{{\mathcal K}}
\def\O{{\mathcal O}}
\def\Q{{\mathcal Q}}
\def\G{{\mathcal G}}

\def\R{I\!\!R}

\def\ov{\overline}
\def\Tilde{\widetilde}
\def\forall{\hbox{for all }~}

\def\v{\vskip 1em}
\def\vs{\vskip 2em}
\def\vsk{\vskip 4em}

\def\be{\begin{equation}}
\def\beq{\begin{equation}}
\def\bel{\begin{equation}\label}
\def\eeq{\end{equation}}

\def\vp{\varphi}

\def\({\left(\begin{array}{cccccc}}
\def\){\end{array}\right)}

\def\({\left(\begin{array}{cccccc}}
\def\){\end{array}\right)}

\def\bes{\begin{eqnarray}}
\def\ees{\end{eqnarray}}

\newcommand{\bea}{\begin{eqnarray}}
\newcommand{\eea}{\end{eqnarray}}
\newcommand{\beann}{\begin{eqnarray*}}
\newcommand{\eeann}{\end{eqnarray*}}

\newcommand{\bp}{\begin{proof}}
\newcommand{\ep}{\end{proof}}

\begin{document}
\title{\bf Unique Conservative Solutions to a Variational Wave Equation}
\vs

\author{Alberto Bressan$^{(*)}$, Geng Chen$^{(**)}$, and Qingtian Zhang$^{(*)}$\\    \\
{\small (*) Department of Mathematics, Penn State University,
University Park, Pa.~16802, U.S.A.}\\
{\small (**) School of Mathematics,
Georgia Institute of Technology,
Atlanta, Ga.~30332, U.S.A.}
\\
\,
\\
{\small e-mails:~ bressan@math.psu.edu~,~gchen73@math.gatech.edu~, ~zhang\_q@math.psu.edu}}

\maketitle

\begin{abstract} Relying on the analysis of characteristics,
we prove the
uniqueness of conservative solutions to the
variational wave equation $u_{tt}-c(u) (c(u)u_x)_x=0$.
Given a  solution $u(t,x)$, even if the wave speed
$c(u)$ is only H\"older continuous in the $t$-$x$ plane, one can still define
forward and backward characteristics in a unique way.
Using a new set of independent variables $X,Y$,
constant along characteristics,
we prove that
$t,x,u$,  together with other variables, satisfy a semilinear system
with smooth coefficients.
From the uniqueness of the solution to this semilinear system,
 one obtains the uniqueness of
conservative solutions
to the Cauchy problem for the wave equation
with general initial data $u(0,\cdot)\in H^1(\R)$,
$u_t(0,\cdot)\in L^2(\R)$.

\end{abstract}
\v
\section{Introduction}
\setcounter{equation}{0}

Consider the Cauchy problem for the quasilinear
second order wave equation
\beq
u_{tt} - c(u)\big(c(u) u_x\big)_x=0\,,\label{1.1}
\eeq
with initial data
\beq
u(0,x)=u_0(x)\,,\qquad
u_t(0,x)=u_1(x)\,.\label{1.2}
\eeq
Here
$u_0\in H^1(\R)$ while $u_1\in \L^2(\R)$.
We assume that the wave speed $c:\R\mapsto \R_+$ is
a smooth, bounded, uniformly positive function, satisfying
\bel{cprop}
0~<~c_0~\leq ~c(u)~<~M,\qquad\qquad   |c'(u)|~<~ M\qquad\quad\forall u.\eeq
The above Cauchy problem  has been studied in several papers \cite{BH, BZ,
GHZ, ZZ03, ZZ05}.
In particular, the analysis in
\cite{BZ} shows that the problem (\ref{1.1})-(\ref{1.2}) has a weak solution
which conserves the total energy. Indeed, a global flow
of such solutions can be constructed, both forward and backward in time,
exhibiting some kind of continuous dependence on the initial data.
The approach developed in \cite{BZ} relied on the introduction of a set of auxiliary
variables.  Using these variables, one obtains a semilinear system of equations
having unique solutions.   In terms of the original variables, this yields a
 solution $u=u(t,x)$ of the Cauchy problem (\ref{1.1})-(\ref{1.2}), for which the
 total energy is a.e.~conserved.   The main results in \cite{BZ} can be
 summarized as follows.
\v
\n{\bf Theorem 1.} {\it Let $c:\R\mapsto \R$ be
a smooth function satisfying (\ref{cprop}).  Assume that the initial data
$u_0$ in (\ref{1.2}) is absolutely continuous, and that
$(u_0)_x\in\L^2\,$,   ~$u_1\in\L^2$.
Then the Cauchy problem (\ref{1.1})-(\ref{1.2}) admits a weak solution
$u=u(t,x)$, defined for all $(t,x)\in\R\times\R$.
In the $t$-$x$ plane, the function $u$ is locally H\"older continuous
with exponent $1/2$.   This solution $t\mapsto u(t,\cdot)$ is
continuously differentiable as a map with values in $\L^p_{\rm loc}$, for all
$1\leq p<2$. Moreover, it is Lipschitz continuous w.r.t.~the $\L^2$
distance, i.e.
\beq
\big\|u(t,\cdot)-u(s,\cdot)\big\|_{\L^2}~\leq~ L\,|t-s|\label{1.lip}
\eeq
for all $t,s\in\R$.
The equation (\ref{1.1}) is satisfied in integral sense, i.e.
\beq
\dint \Big[\phi_t\, u_t - \big(c(u) \phi\big)_x c(u)\,u_x\Big]\,dxdt ~=~0
\label{consol}
\eeq
for all test functions $\phi\in\C^1_c$.
Moreover, the maps $t\mapsto u_t(t, \cdot)$ and $t\mapsto u_x(t,\cdot)$ are continuous
with values in $\L^p_{loc}(\R)$, for every $p\in [1,2[\,$.}
\v
In general, the solution constructed in Theorem 1 is not unique.
To  select a unique solution, additional properties
must be imposed. In particular, one can require that the total energy be conserved.

It is convenient to introduce the variables
\beq
\left\{
\begin{array}{rcl}
R & \doteq  &u_t+c(u)u_x\,, \\
S & \doteq  &u_t-c(u)u_x\,,
\end{array} \right.\label{2.1}
\eeq
so that
\beq
u_t={R+S\over 2}\,,\qquad\qquad u_x={R-S\over 2c}\,.\label{2.2}
\eeq
By (\ref{1.1}), the variables $R,S$ satisfy
\beq
\left\{
\begin{array}{rcl}
R_t-cR_x &= & {c'\over 4c}(R^2-S^2), \\ [3mm]
S_t+cS_x &= & {c'\over 4c}(S^2-R^2).
\end{array} \right.
\label{2.3}
\eeq
Multiplying the first equation in (\ref{2.3}) by $R$ and the
second one by $S$, one obtains balance laws for $R^2$ and $S^2$, namely
\beq\left\{
\begin{array}{rcl}
(R^2)_t - (cR^2)_x & = & {c'\over 2c}(R^2S - RS^2)\, , \\ [3mm]
(S^2)_t + (cS^2)_x & = & - {c'\over 2c}(R^2S -RS^2)\,.
\end{array}
\right.\label{2.4}
\eeq
As a consequence, the following quantities are conserved:
\beq
E~\doteq~ 2\big(u_t^2+c^2u_x^2\big)~=~R^2+S^2\,,\qquad\qquad
\mathcal{M}~\doteq ~-u_tu_x~=~{S^2-R^2\over 4c}\,.\label{2.5}
\eeq
One can think of $R^2$ and $S^2$ as the energy of backward and
forward moving waves, respectively.  Notice that these are not separately conserved.
Indeed, by (\ref{2.4}) energy is transferred from forward to backward waves, and
vice versa.

{\bf Theorem 2.} {\it Under the previous assumptions, a solution $u=u(t,x)$
can be constructed which is {\em conservative} in the following sense.

 There exists two  families of positive Radon measures
 on the real line: $\{\mu_-^t\}$ and $\{\mu_+^t\}$, depending continuously
on $t$ in the weak topology of measures, with the following properties.
\begi
\item[(i)] At every time $t$ one has
\bel{E0}\mu_-^t(\R)+\mu_+^t(\R)~=~E_0~\doteq~2
\int_{-\infty}^\infty \Big[u_1^2(x) + \bigl(c(u_0(x)) u_{0,x}(x)\bigr)^2\Big]\, dx \,.\eeq

\item[(ii)] For each $t$, the absolutely continuous parts of $\mu_-^t$ and
$\mu_+^t$ w.r.t.~the Lebesgue measure
have densities  respectively given  by
\bel{dmu}R^2 ~=~
\bigl(u_t + c(u) u_x\bigr)^2,\qquad\qquad
S^2~=~\bigl(u_t - c(u) u_x\bigr)^2.
\eeq
\item[(iii)]  For almost every $t\in\R$, the singular parts of $\mu^t_-$ and $\mu^t_+$
are concentrated on the set where $c'(u)=0$.
\item[(iv)] The measures $\mu_-^t$ and $\mu_+^t$ provide measure-valued solutions
respectively to the balance laws
\bel{mb4}\left\{
\begin{array}{rcl}
w_t - (cw)_x & = & {c'\over 2c}(R^2S - RS^2)\, , \\ [3mm]
z_t + (cz)_x & = & - {c'\over 2c}(R^2S -RS^2)\,.
\end{array}
\right.
\eeq
\endi
}

{\bf Remark 1.} In principle, the equations (\ref{mb4}) should be written as
\bel{mb5}\left\{
\begin{array}{rcl}
w_t - (cw)_x & = & {c'\over 2c}(Sw - Rz)\, , \\ [3mm]
z_t + (cz)_x & = & - {c'\over 2c}(Sw -Rz)\,.
\end{array}
\right.
\eeq
This reflects the fact that, if $w= w^a + w^s$ is a measure with an absolutely continuous
and a singular part, then both of these multiply $S$.
However, we are here
making the assumption that the solution is {\it conservative}, so that by (iii)
the product $c'(u) w^s =0$
for a.e.~time $t$.     For this reason, on the right hand side
of (\ref{mb5}) we can  replace $w$ with the measure $w^a$
having density $R^2$ w.r.t.~Lebesgue measure.
Similarly, we can  replace $z$ with the measure $z^a$
having density $S^2$ w.r.t.~Lebesgue measure.
\v
Observe that the total energy represented by the sum $\mu_-^t+\mu_+^t$ is
conserved in time.  Occasionally, some of this energy is concentrated
on a set of measure zero. At the times $\tau$ when this happens,
$\mu_-^t+\mu_+^t$ has a non-trivial singular part and
$$E(\tau)~\doteq~\int_{-\infty}^\infty\Big[u_t^2(\tau,x) +
c^2(u(\tau,x)) u^2_x(\tau,x)
\Big]\, dx  ~<~E_0\,.$$
The condition (iii) puts some restrictions
on the set of such times $\tau$. In particular, if $c'(u)\not= 0$
for all $u$, then this set has measure zero.

Our present goal is to understand whether these conservative solutions are unique.
In a way, our approach is
the inverse of \cite{BZ}.
Given a conservative solution
$u=u(t,x)$, we define a set of independent variables $X,Y$ and dependent variables
$u,w,z,p,q$, and show that these satisfy a suitable semilinear system of equations.
By proving that this semilinear system has unique solutions, we eventually
obtain the uniqueness of solutions to the original equation (\ref{1.1}).

In essence, this semilinear system describes the evolution of $u$ and its derivatives
along characteristic
curves, i.e.~curves $t\mapsto x^\pm(t)$ which satisfy the ODEs
\bel{dxpm}\dot x^-(t)~=~-c(u(t, x^-(t))),\qquad\qquad \dot x^+(t)~=~c(u(t, x^+(t))).
\eeq
At this naive level, the approach runs into a fundamental
difficulty.  Namely,  since the solution $u$ is only H\"older continuous, for a given $\bar y\in \R$
the Cauchy problems for the ODEs in (\ref{dxpm}) with initial data
\bel{icy}
x^\pm(0)~=~\bar y,\eeq
may well have multiple solutions.
To overcome this difficulty, our analysis relies on two ideas.  To simplify the exposition,
we here assume that the measures $\mu^t_-, \mu^t_+$ are absolutely
continuous for a.e.~$t$.
\begi
\item{}  The two balance laws (\ref{2.4}) imply
\bel{ceq1}
{d\over dt}  \int_{-\infty}^{x^-(t)}R^2(t,x)\, dx~=~\int_{-\infty}^{x^-(t)} {c'\over 2c} (R^2S-RS^2)
\, dx\,,\eeq
\bel{ceq2}{d\over dt}  \int_{-\infty}^{x^+(t)}S^2(t,x)\, dx~=~-\int_{-\infty}^{x^+(t)} {c'\over 2c} (R^2S-RS^2)
\, dx\,.\eeq
While the Cauchy problem (\ref{dxpm})-(\ref{icy}) can have multiple solutions,
the characteristic curves $t\mapsto x^\pm(t)$ can be uniquely determined
by combining all the
equations in (\ref{dxpm})--(\ref{ceq2}).

\item{}  Instead of the variables $(t,x)$, it is convenient to work with
an adapted set of variables $x(t, \alpha)$, $y(t,\beta)$, where
\bel{xbeta} x(t,\alpha) + \int_{-\infty}^{x(t,\alpha)} R^2(t,\xi)\, d\xi~=~\alpha\,,\eeq
\bel{ybeta} y(t,\beta) + \int_{-\infty}^{y(t,\beta)} S^2(t,\xi)\, d\xi~=~\beta\,.\eeq
Here the parameter $\alpha$ singles out  a  backward characteristic, while $\beta$
singles out a forward characteristic.
\endi

Our main result is the following.

{\bf Theorem 3. } {\it Let $c:\R\mapsto \R$ be
a smooth function satisfying (\ref{cprop}).  For any initial data
$u_0\in H^1(\R)\,$,   ~$u_1\in\L^2(\R)$,
the conservative solution to  Cauchy problem (\ref{1.1})-(\ref{1.2}),
which satisfies all conditions
(i)--(iv) in Theorem 2, is unique.}
\v
The main technique used in the proof is similar to the paper
\cite{BCZ} by the same authors, dealing with the
Camassa-Holm equation, and was inspired by
the uniqueness result in \cite{Daf}.
However, in the present case of a second order
wave equation, the analysis is harder.
Indeed, for the Camassa-Holm equation one has a single family of characteristics.
After a change of variables, each characteristic is obtained by solving an ODE with Lipschitz continuous right hand side.

On the other hand, for the wave equation (\ref{1.1}) one has two families of characteristics moving forward and backward, respectively.
After a change of variables, the ODEs which determine these
characteristics are still no better than H\"older continuous.
However, the singularities are transversal.   Uniqueness and Lipschitz continuous
dependence of solutions on the initial data can thus be established using ideas from
\cite{B1,  BC, BS}.

The paper is organized as follows. In Section~2 we review the basic
equations and prove an a priori estimate on the total amount of wave interactions.
In Section~3 we show that, for a given conservative solution $u=u(t,x)$, one can uniquely
determine a forward and a backward characteristic through each initial point.
In Section~4 we introduce the characteristic coordinates $(X,Y)$ and prove
Lipschitz continuity of map $(X,Y)\mapsto (t,x,u)$. In Section 5 we introduce
some additional variables and show that they satisfy a semilinear system
with smooth coefficients.  The uniqueness of solutions to this semilinear system
yields the uniqueness of conservative solutions to the original wave equation (\ref{1.1}).
\v
\section{Adapted variables and wave interaction estimate}
\setcounter{equation}{0}

Recalling (\ref{cprop}), for notation convenience, we introduce  the constant
\bel{C0}
C_0~\doteq~\left\|\frac{c'}{2c}\right\|_{\L^\infty}~\leq~ \frac{M}{2c_0}\,.
\eeq
Let $u=u(t,x)$ be a conservative solution of (\ref{1.1}),
having all the properties listed in Theorems~1 and 2.
For any time $t$ and  any $\alpha,\beta\in \R$, we define the points $x(t,\alpha)$ and
$ y(t,\beta)$ by setting
\bel{xadef}
x(t,\alpha)~\doteq~\sup\Big\{ x\,;~~x+\mu^t_-\bigl(\,]-\infty,x]\,
\bigr) ~<~\alpha\Big\}\,,\eeq
\bel{ybdef}y(t,\beta)~\doteq~\sup\Big\{ x\,;~~x+\mu^t_+\bigl(\,]-\infty,x]\,\bigr)
 ~<~\beta\Big\}\,.\eeq
Notice that the above holds if and only if, for some $\theta,\theta'\in [0,1]$, one has
\bel{xa}
x(t,\alpha)+\mu_-^t\Big(\bigl]-\infty\,, ~x(t,\alpha)
\bigr[\Big)+\theta\cdot \mu^t_- \Big(\bigl\{x(t,\alpha)
\bigr\}\Big)~=~\alpha\,,\eeq
\bel{yb}
y(t,\beta)+\mu_+^t\Big(\bigl]-\infty\,,~ y(t,\beta)
\bigr[\Big) +\theta'\cdot \mu^t_- \Big(\bigl\{y(t,\beta)
\bigr\}\Big)~=~
\beta\,.\eeq
Since the measures $\mu^t_-$, $\mu^t_+$ are both positive and bounded,
it is clear that these points
are well defined.   In the absolutely continuous case,
the equations (\ref{xa})-(\ref{yb}) are equivalent to (\ref{xbeta})-(\ref{ybeta}).
\v
{\bf Lemma 1.}  {\it For every fixed $t$,  the maps $\alpha\mapsto x(t,\alpha)$  and
 $\beta\mapsto y(t,\beta)$  are both Lipschitz continuous with constant 1.
 Moreover, for fixed $\alpha,\beta$, the maps $t\mapsto  x(t,\alpha)$ and
 $t\mapsto  y(t,\beta)$
are absolutely continuous and locally H\"older continuous with exponent $1/2$.}
\v
{\bf Proof.}
{\bf 1.} The first part is straightforward.  Indeed, if
$$x_1~\doteq~ x(t,\alpha_1)~<~x(t,\alpha_2)~\doteq~x_2\,,$$
then
$$x_2-x_1~\leq~x_2-x_1 +\mu_-^t\bigl(\,]x_1, x_2[\,\bigr)~\leq~\alpha_2-\alpha_1\,.$$
The same argument applies to the map $\beta\mapsto y(t,\beta)$.
\v
{\bf 2.} To prove the second statement, denote by $\mu_-^t\otimes\mu_+^t$
 the product measure on $\R^2$ and
consider the wave interaction potential
\bel{wip}
Q(t)~\doteq~(\mu_-^t\otimes\mu_+^t)\Big( \{(x,y)\,;~~x>y\}\Big).\eeq
We recall that $\mu_-^t(\R)+\mu_+^t(\R)=E_0$ is the total energy, constant in time.
Since $R^2(t,\cdot)$ and $S^2(t,\cdot)$ provide the
absolutely continuous parts of $\mu_-^t$ and $\mu_+^t$, respectively,
recalling (\ref{cprop}) and using the balance laws  (\ref{mb4}) we obtain
\bel{DQ} \bega{rl}\ds{d\over dt}Q(t)&\ds\leq~-c_0\int_{-\infty}^{+\infty}
 S^2(t,x)R^2(t,x)\, dx +
\frac{M}{2c_0}\int\left( \int_x^{+\infty}|R^2 S-RS^2|(t,y)\, dy\right)d\mu_-^t(x) \cr\cr
&\qquad\qquad \ds+
\frac{M}{2c_0}\int \left( \int_{-\infty}^y|R^2 S-RS^2|(t,x)\, dx\right)d\mu_+^t(y)\cr\cr
&\ds\leq~-c_0\int_{-\infty}^{+\infty} S^2R^2\, dx +
\frac{M}{2c_0} \bigr(\mu_-^t(\R)+\mu_+^t(\R)\bigl) \int_{-\infty}^{+\infty}|R^2 S-RS^2|\, dx\cr\cr
&\ds\leq~-c_0\int_{-\infty}^{+\infty}  S^2R^2\, dx +
\frac{M}{2c_0} E_0 \int_{-\infty}^{+\infty}\bigl(|R^2 S| +|RS^2|\bigr)\, dx\cr\cr
&\ds\leq~-{c_0\over 2}\int_{-\infty}^{+\infty}  S^2R^2\, dx +
 \int_{\{2ME_0>c_0^2|S|\}}|R^2 S|\, dx +
\int_{\{2ME_0>c^2_0|R|\}}|RS^2|\, dx\cr\cr
&\ds\leq ~ -{c_0\over 2}\int_{-\infty}^{+\infty}  S^2R^2\, dx +
{2ME_0\over c_0^2} \int_{-\infty}^{+\infty}(R^2 +S^2)\, dx \cr\cr
&\ds
\leq~-{c_0\over 2}\int_{-\infty}^{+\infty}  S^2R^2\, dx +
{2ME_0^2\over c^2_0}\,.
\enda
\eeq
Since $Q(t)\leq E_0^2$ for every time $t$, from (\ref{DQ})
it follows
\bel{SRb}\int_0^T \int_{-\infty}^{+\infty} R^2(t,x) S^2(t,x)\, dx\, dt~\leq~
{2\over c_0}\cdot\left[ \bigl( Q(0)-Q(T)\bigr) + {2ME_0^2\over c^2_0}\, T\right]~\leq~
\frac{2E_0^2}{c_0}+\frac{4ME_0^2}{c_0^3}\,T\,.\eeq

\v
{\bf 3.} For a given $\tau$ and any $\ve\in\,]0,1]$, we now estimate
\bel{dxy} \bega{l}\ds
\int_\tau^{\tau+\ve} \int_{-\infty}^{+\infty} |R^2 S-RS^2|
\,dx\, dt\cr
\cr
\qquad \ds\leq~ \int_\tau^{\tau+\ve}\int_{S\le\ve^{-1/2}}|R^2 S|\,dx\, dt+
\int_\tau^{\tau+\ve}\int_{R\le\ve^{-1/2}}|R S^2|\,dx\, dt\cr\cr
\ds\qquad \qquad  + \int_\tau^{\tau+\ve}\int_{S>\ve^{-1/2}}|R^2 S|\,dx\, dt
 +\int_\tau^{\tau+\ve}\int_{R>\ve^{-1/2}}|R S^2|\,dx\, dt\cr\cr
\ds\qquad
 \leq~ \ve^{-1/2} \int_\tau^{\tau+\ve}\int_{-\infty}^{+\infty} (R^2+S^2) \,dx\, dt
\cr\cr
\ds\qquad\qquad
 + \int_\tau^{\tau+\ve}\int_{S\geq\ve^{-1/2}}(R^2S^2) {|R^2 S|\over R^2 S^2}\,dx\, dt
 + \int_\tau^{\tau+\ve}\int_{R\geq\ve^{-1/2}}(R^2S^2) {|R S^2|\over R^2 S^2}\,dx\, dt
\cr\cr
\ds\qquad\leq~
E_0 \ve^{1/2} + 2\ve^{1/2} (\frac{2E_0^2}{c_0}+\frac{4ME_0^2}{c_0^3}\ve)\,\cr\cr
\ds\qquad\leq~\left[E_0+ \frac{4E_0^2}{c_0}+\frac{8ME_0^2}{c_0^3}
\right]\ve^{1/2} .
 \enda
 \eeq
 As a consequence, the function $\zeta$ defined by
 \bel{zdef}
 \zeta(\tau)~\doteq~\int_0^\tau \int_{-\infty}^{+\infty}
 \left|\frac{c'}{2c}(R^2S-RS^2)\right|
\,dx\, dt\eeq
is locally H\"older continuous, nondecreasing, with sub-linear growth.
Since $R^2S-RS^2\in \L^1([0,T]\times \R)$, by Fubini's theorem the map
$t\mapsto\ds \int
 \left|\frac{c'}{2c}(R^2S-RS^2)\right|
\,dx$ is  in $\L^1([0,T])$.   By its definition at (\ref{zdef}),
the function $\zeta$ is absolutely continuous.
Recalling (\ref{cprop}), for $0<t_2-t_1\leq 1$ we have
$$
\zeta(t_2)-\zeta(t_1)~\leq ~C_1(t_1-t_2)^{1/2},
$$
where the constant $C_1$ is defined as
\bel{const1}
C_1~\doteq~ {M\over c_0}\, \left[E_0+ \frac{4E_0^2}{c_0}+\frac{8ME_0^2}{c_0^3}
\right].
\eeq
\v
{\bf 4.} We recall that the family of measures
$\mu_-^t$ satisfies the balance law in (\ref{mb4}) with velocity $-c(u)\in [-M,0]$.
For any $t_1<t_2$ and any $\alpha$, this yields the inequalities
\bel{i1}
\mu_-^{t_2}\Big( \bigl]-\infty\,,~ x(t_1,\alpha)\bigr[\Big)
~\geq~\mu_-^{t_1}\Big( \bigl]-\infty\,,~ x(t_1,\alpha)\bigr[\Big) -
\bigl[\zeta(t_2)
-\zeta(t_1)\bigr]\,,\eeq
\bel{i2}
\mu_-^{t_2}\Big( \bigl]-\infty\,,~ x(t_1,\alpha)-M(t_2-t_1)\bigr[\Big)
~\leq~\mu_-^{t_1}\Big( \bigl]-\infty\,,~ x(t_1,\alpha)\bigr[\Big) +\bigl[\zeta(t_2)
-\zeta(t_1)\bigr]\,.\eeq
{} From the definition (\ref{xa}) it thus follows
\bel{x11}
x(t_1,\alpha)-M(t_2-t_1) -  \bigl[\zeta(t_2)
-\zeta(t_1)\bigr]  ~\leq~ x(t_2,\alpha) ~\leq ~x(t_1,\alpha)+\bigl[\zeta(t_2)
-\zeta(t_1)\bigr].\eeq
By the properties of the function $\zeta$, proved in step {\bf 3}, this achieves the proof.
Of course, the same argument can be applied to the map $t\mapsto y(t,\beta)$.
\endproof

{\bf Remark 2.}   For each fixed $t$,  the map $\alpha\mapsto x(t,\alpha)$ is Lipschitz continuous, hence a.e.~differentiable.   We can define the set of
singular points $\Omega^t$ and the set of singular values $V^t$ according to
\bel{sst}\Omega^t~\doteq~\Big\{ \alpha\in\R\,;~~
{\partial\over\partial\alpha} x(t,\alpha)~=~0~~\hbox{or else this partial derivative does not exist}\Big\},\eeq
\bel{Vvt}V^t~\doteq~\Big\{ x(t,\alpha)\,;~~\alpha\in \Omega^t\Big\}.\eeq
Observe that $V^t$ has zero Lebesgue measure.

In general, the map $\alpha\mapsto x(t,\alpha)$ is onto but not one-to-one.
However, for every regular value $x_0\in \R\setminus V^t$
there exists a unique $\alpha_0$ such that $x_0 = x(t,\alpha_0)$.

If now $f\in\L^1(\R)$,
the composition
$\tilde f(\alpha)~\doteq~f(x(t,\alpha))$ is well defined
for a.e.~$\alpha\in\R\setminus \Omega^t\,$.  The integral
of $f$ can be computed by a change of variables:
\bel{cvi}
\int f(x)\, dx~=~\int_{\R\setminus\Omega^t} \tilde f(\alpha) \cdot
{\partial\over\partial\alpha} x(t,\alpha)\, d\alpha~=~
\int_{\R\setminus\Omega^t} \tilde f(\alpha) \cdot
{1\over 1+R^2(t, x(t,\alpha))}\, d\alpha\,.\eeq

\section{Recovering the characteristic curves}
\setcounter{equation}{0}

The next lemma, which plays a crucial role in our analysis,
shows that for a conservative solution
the characteristic curves can be uniquely determined.
Observe that, in the general case where the measures $\mu^t_-,\mu^t_+$
need not be absolutely continuous, the identities (\ref{ceq1})-(\ref{ceq2})
can be written in the equivalent integrated form
\bel{be1}
\int_{-\infty}^{\bar y} R^2(0,x)\, dx + \int_0^t
\int_{-\infty}^{x^-(s)} {c' (R^2S-RS^2)\over 2c}
\, dx\, ds~=~\mu^t_-\Big( ]-\infty, x^-(t)[\Big) + \theta(t,\bar y)\cdot
\mu^t_-\Big( \{ x^-(t)\}\Big),\eeq
\bel{be2}
\int_{-\infty}^{\bar y} S^2(0,x)\, dx - \int_0^t
\int_{-\infty}^{x^+(s)} {c' (R^2S-RS^2)\over 2c}
\, dx\, ds~=~\mu^t_+\Big( ]-\infty, x^-(t)[\Big) + \theta'(t,\bar y)\cdot
\mu^t_+\Big( \{ x^+(t)\}\Big),\eeq
for some functions $\theta, \theta'\in [0,1]$.
\v
{\bf Lemma 2.} {\it Let $u$ be a conservative solution of (\ref{1.1}), satisfying
the properties stated in Theorems 1 and 2.
Then, for any $\bar y\in \R$, there exists unique Lipschitz continuous maps
$t\mapsto x^\pm(t)$ which satisfy (\ref{dxpm})-(\ref{icy})
together with
(\ref{be1})-(\ref{be2}).}
\v
{\bf Proof.}  We claim  that there exists a unique function
$t\mapsto \alpha(t)$ such that
\bel{xal}x^-(t)~=~x(t,\alpha(t))\eeq
satisfies the equations in (\ref{dxpm})-(\ref{icy}) and (\ref{ceq1}).
It suffices to prove the claim on the time interval $t\in [0,1]$, then iterate the argument
by induction.
The proof will be given in several steps.
\v
{\bf 1.}
Integrating the first equation in (\ref{dxpm}) w.r.t.~time and summing it with
(\ref{be1})
we obtain
\bel{ceq3}\bega{l}\ds
x^-(t)+\mu_-^{t}\Big( \bigl]-\infty\,,~ x^-(t)\bigr[\Big)+
\theta(t) \cdot\mu_-^{t}\Big( \{x^-(t)\}\Big)
\cr\cr\ds\qquad =~\bar y + \int_{-\infty}^{\bar y}R^2(0,x)\, dx
+ \int_0^t\left( -c(u(s, x^-(s)))+\int_{-\infty}^{x^-(s)} {c' (R^2S-RS^2)\over 2c}
\, dx\right) ds\,,\enda\eeq
for some $\theta(t)\in [0,1]$.

From
(\ref{xal}) and (\ref{ceq3})
we  obtain an integral equation for $\alpha$,
namely
\bel{ceq4}
\alpha(t)~=~\bar\alpha + \int_0^t \left(-c(u(s,x^-(s)))+\int_{-\infty}^{x(s, \alpha(s))}
{c' (R^2S-RS^2)\over 2c}
\, dx\right) ds\,.\eeq
Here
\bel{ica}\bar \alpha~=~\alpha(0)~=~\bar y +\int_{-\infty}^{\bar y} R^2(0,x)\, dx\,.\eeq
Notice that the equation (\ref{ceq4}) is equivalent to
\bel{ceq5}\bega{rl}
\dot \alpha(t)&\ds=~G(t, \alpha(t))~\doteq~  -c\bigl(u(t, x(t,\alpha(t)))\bigr)+\int_{-\infty}^{x(t, \alpha(t))}
{c' (R^2S-RS^2)\over 2c}
\,dx\,,\enda\eeq
with initial data (\ref{ica}).  We take (\ref{ceq4}) as the starting point for our analysis.
In the following steps we will show that this integral equation has a unique solution
$t\mapsto \alpha(t)$.
Moreover, the function $t\mapsto x^-(t) = x(t,\alpha(t))$  satisfies the first equation in (\ref{dxpm}), as well as
(\ref{ceq1}).
\v
{\bf 2.} We first prove the existence of a solution to
(\ref{ceq4}) on the time interval $[0,1]$.
Consider the Picard map $\P: \C^0([0,1])\mapsto \C^0([0,1])$, defined as
\bel{Pdef}
\P\alpha(t)~\doteq~\bar\alpha +\int_0^t\left[-c\bigl(u(t, x(t,\alpha(t)))\bigr)+\int_{-\infty}^{x(t, \alpha(t))}
{c' (R^2S-RS^2)\over 2c}
\,dx\,\right] ds\,.
\eeq
We claim that $\P$ is a continuous transformation of
a compact convex  set $\K\subset \C^0([0,1])$ into itself, with the usual norm
$$\|f\|_{\C^0}~\doteq ~\max_{t\in [0,1]}|f(t)|.$$
 Here  the set $\K$ is a set of H\"older continuous functions, defined by
\bel{Kdef}
\ds \K~\doteq~\{f\in \C^{1/2}([0,1])\,;~~~\|f\|_{\C^{1/2}}~\leq~ C_K,
~~~ f(0)=\bar\alpha\},
\eeq
for a suitable constant $C_K$, to be determined later.
Indeed, for any $t\in[0,1]$ one has
\bel{in8}\bega{ll}\ds
&|\P\alpha_1(t)-\P\alpha_2(t)|\cr\cr
=& \ds\left|\int_0^t \bigl(c(u(s,x(s,\alpha_1(s))))-c(u(s,x(s,\alpha_2(s))))\bigr)+ \int_{x(s,\alpha_1(s))}^{x(s,\alpha_2(s))}{c' (R^2S-RS^2)\over 2c}\, dxds\right|\,\cr\cr
\leq &\ds
\int_0^t M \|u\|_{C^{1/2}}|\alpha_2-\alpha_1|^{1/2} ds
+\int_0^t\int_{x(s,\alpha_1(s))}^{x(s,\alpha_2(s))}{|c'|\over 2c}\Big(|RS|
\bigl(|R|+|S|\bigr)\Big)\,dxds\,\cr\cr
\leq &~\ds C_0\int_0^1 E_0^{1/2}|\alpha_2-\alpha_1|^{1/2}ds\cr\cr
&\ds~\qquad +C_0\left(\int_0^1\int_{x(s,\alpha_1(s))}^{x(s,\alpha_2(s))}R^2S^2dxds\right)^{1/2}
\left(\int_0^1\int_{x(s,\alpha_1(s))}^{x(s,\alpha_2(s))}(R^2+S^2)dxds\right)^{1/2}.
\enda
\eeq
Here we have used H\"older's inequality and the fact that $\alpha\mapsto x(t,\alpha)$ is Lipschitz continuous with constant 1.

Consider any function
$\alpha_1\in \C^0([0,1])$.
As $\|\alpha_2-\alpha_1\|_{\C^0}\to 0$ we have $\|x(\cdot ,\alpha_2)-x(\cdot,\alpha_1)\|_{\C^0}\to 0$ as well.  Since
the functions $R^2S^2$, $R^2$, and $S^2$ are all in $\L^1([0,1]\times\R)$,
the right hand side of (\ref{in8}) approaches zero.  This proves the continuity
of the map $\P$, in the $\C^0$ norm.

Next, we need to show that, for a suitable choice of $C_K$, the
transformation  $\P$ maps the compact convex set $\K$ in (\ref{Kdef}) into itself.
For $0<t_2-t_1\leq 1$, recalling (\ref{zdef}) we obtain
\bel{Hcont}\bega{rl}
\bigl|\P\alpha(t_2)- \P\alpha(t_1)\bigr|&\ds=~\left|\int_{t_1}
^{t_2}\left(-c(u(s,x(s,\alpha(s))))+\int_{-\infty}^{x(s,\alpha(s))}
{c' (R^2S-RS^2)\over 2c}\, dx\right) ds\right|\cr\cr
&\ds \leq~ \int_{t_1}^{t_2} M \, dt+
\int_{t_1}^{t_2}\int_{-\infty}^\infty\left|\frac{c'}{2c}(R^2S-RS^2)\right|dx dt\cr\cr
&\leq~\ds M\cdot (t_2-t_1) + [\zeta(t_2)-\zeta(t_1)]\cr\cr
&\leq~\ds (M+C_1)(t_2-t_1)^{1/2}.
\enda
\eeq
The above computation also yields
\bel{Pai}\max_{t\in [0,1]}~ |\P\alpha(t)|~\leq~ \bar \alpha +M+C_1\,.\eeq
Together, (\ref{Hcont}) and (\ref{Pai})  yield
an a priori bound on the H\"older norm
$$\|\P\alpha\|_{\C^{1/2}}~\leq~\bar \alpha +2(M+C_1),$$
where  $C_1$ is the constant in \eqref{const1}.

By (\ref{Pdef}),  $\P\alpha(0)=\bar\alpha$.
Choosing the constant $C_K\doteq \bar \alpha +2(M+C_1)$
in (\ref{Kdef}), we obtain  that $\P$ maps $\K$ into $\K$.

By Schauder's fixed point theorem,  the integral equation  (\ref{ceq4})
has at least one solution.  Iterating the argument, this solution can be
extended to any time interval $t\in [0,T]$.
\v
{\bf 3.} In this step and the next one we prove that  $x^-(\tau)\doteq x(\tau,\alpha(\tau))$
satisfies the first equation in (\ref{dxpm}) at a.e. time $\tau$.

Since $\frac{c'}{2c}[R^2S-RS^2]\in \L^1([0,T]\times\R)$,
a classical theorem of Lebesgue implies that
$$
\lim_{r\to 0+}~{1\over \pi r^2}\,\dint_{\{(\tau-t)^2+(y-x)^2\leq r^2\}}
\frac{c'}{2c}[R^2S-RS^2](\tau,y)dyd\tau~=~\frac{c'}{2c}[R^2S-RS^2](t,x),
$$
for all $(t,x)\in ~]0,T[\,\times \R$ outside a null
set $\N_2$ whose 2-dimensional measure is zero.

If one divides by $r$ instead of $r^2$, by Corollary 3.2.3 in  \cite{Zi}
there is a set $\N_1\subset \N_2$
whose 1-dimensional Hausdorff measure is zero and such that$$
\limsup_{r\to 0+}~{1\over r}\,\dint_{\{(\tau-t)^2+(y-x)^2\leq r^2\}}
\frac{c'}{2c}[R^2S-RS^2](\tau,y)dyd\tau~=~0
$$
for every $(t,x)\not\in \N_1$.
Therefore, there exists a null 1-dimensional set $\N\subset [0,T]$
with the properties
\begi
\item[(i)]  For every $\tau\notin\N$ and $x\in \R$ one has $(\tau,x)\notin \N_1$\,,
\item[(ii)] If $\tau\notin \N$ then the map $t\mapsto \zeta(t)$ in (\ref{zdef})
is differentiable at $t=\tau$. Moreover, $\tau$
is a Lebesgue point of the derivative $\zeta'$.
\endi
We conclude this step by observing that the map $t\mapsto
x^-(t)\doteq x(t,\alpha(t))$ is absolutely continuous.  Indeed,
consider a finite sequence of times such that
$$0<s_1<t_1<s_2<t_2<\cdots < s_\nu< t_\nu\,.$$
Using (\ref{x11}) and the fact that the map $\alpha\mapsto x(t,\alpha)$ is contractive,
we obtain
\bel{aca}\bega{rl}S&\ds\doteq~\sum_{k=1}^\nu |x^-(t_k) - x^-(s_k)| ~
=~\sum_{k=1}^\nu \bigl|x^-(t_k, \alpha(t_k)) - x^-(s_k,
\alpha(s_k))\bigr|\cr\cr
&\ds\leq~ \sum_{k=1}^\nu
\bigl|x^-(t_k, \alpha(t_k)) - x^-(t_k,
\alpha(s_k))\bigr|+ \sum_{k=1}^\nu
\bigl|x^-(t_k, \alpha(s_k)) - x^-(s_k,
\alpha(s_k))\bigr|\cr\cr
&\ds\leq~ \sum_{k=1}^\nu
\bigl|\alpha(t_k) -
\alpha(s_k))\bigr|+ \sum_{k=1}^\nu \Big(M(t_k-s_k)+
\bigl|\zeta(t_k) - \zeta(s_k)\bigr|\Big)\,.\enda
\eeq
Since the two maps $\alpha(\cdot)$ and $\zeta(\cdot)$ are both absolutely continuous,
given $\ve>0$ there exists $\delta>0$ such that the inequality
$$ \sum_{k=1}^\nu |t_k-s_k|~\leq~\delta$$
implies that both summations on the right hand side of (\ref{aca}) are $<\ve/2$.
This proves the absolute continuity of the map $t\mapsto x^-(t)$.

By possibly enlarging the null set $\N\subset [0,T]$ we can assume that,
in addition to (i)-(ii) above, one has
\begi
\item[(iii)] The functions $t\mapsto x^-(t)$ and
$t\mapsto \alpha(t)$ are differentiable at each point   $\tau\in [0,T]\setminus\N$.
Moreover, each point $\tau \notin \N$ is a Lebesgue point of the derivatives $\dot x^-$ and $\dot\alpha$.
\endi

\begin{figure}[htbp]
\centering
\includegraphics[width=0.7\textwidth]{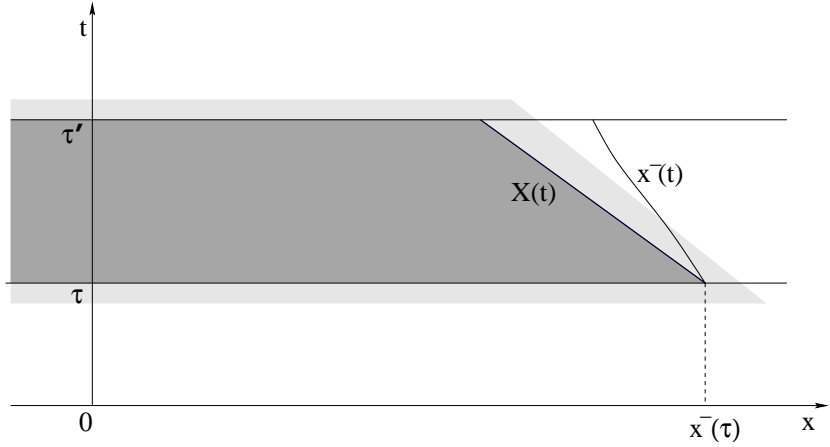}
\caption{The construction used to prove that $\dot x^-(\tau) = -c(u(\tau, x^-(\tau))$.
The shaded area is the support of the test function $\vp^\epsilon$. }
\label{f:wa28}
\end{figure}

{\bf 4.} Let now $\tau\not\in \N$.
 We claim  that the map $t\mapsto x^-(t)\doteq x(t,\alpha(t))$
satisfies the first equation in (\ref{dxpm}) at time $t=\tau$.

Assume, on the contrary, that $\dot x^-(\tau)\not= -c(u(\tau, x^-(\tau)))$.
To fix the ideas, let
 \bel{ass}\dot x^-(\tau)~=~ -c(u(\tau, x^-(\tau)))+2\ve_0\eeq
 for some $\ve_0>0$.   The case $\ve_0 <0$ is entirely similar.
To derive a contradiction we choose $\delta>0$ small enough so that, as shown in
Fig.~\ref{f:wa28},
\bel{xpmt}X(t)~\doteq~
x^-(\tau) +(t-\tau) [-c(u(\tau, x^-(\tau)))+\ve_0 ]~<~x^-(t)\eeq
for all $t\in ]\tau, \tau+\delta]$.
Since $u$ is continuous while $u_t,u_x\in\L^2$, by an
approximation argument the identity in (\ref{consol}) remains
valid for any Lipschitz continuous function $\vp$ with compact support.
Given $\tau<\tau'<\tau+\delta$, for
$\epsilon>0$ small we shall construct a Lipschitz approximation $\vp^\epsilon$ to
the characteristic function of the set
$$\Omega~\doteq~\{(t,y)\,;~~t\in [\tau,\tau']\,,~~~y\in [\epsilon^{-1}, ~X(t)]\bigr\}.$$
For this purpose, consider the functions
$$\rho^{\epsilon}(t,y)~\doteq~\left\{\bega{cl} 0 \qquad &\hbox{if}\quad y
\leq -\epsilon^{-1},\cr
\epsilon^{-1}(y+\epsilon^{-1}) \qquad &\hbox{if}\quad -\epsilon^{-1}\leq y\leq
\epsilon-\epsilon^{-1},\cr
1\qquad &\hbox{if}\quad \epsilon-\epsilon^{-1} \leq y\leq X(t),\cr
1-\epsilon^{-1}(y-X(t))\qquad &\hbox{if}\quad  X(t)\leq y\leq X(t)+\epsilon,\cr
0\qquad &\hbox{if}\quad y\geq X(t)+\epsilon,\enda\right.$$
\bel{timtest}
\chi^\epsilon(t)~\doteq~\left\{\bega{cl} 0\qquad
&\hbox{if}\quad t\leq \tau-\epsilon,\cr
\epsilon^{-1}(s-\tau+\epsilon)\qquad &\hbox{if}\quad \tau-\epsilon\leq t\leq \tau,\cr
1\qquad &\hbox{if}\quad \tau\leq t\leq  \tau',\cr
1-\epsilon^{-1}(t-\tau') \qquad &\hbox{if}\quad \tau'\leq t<\tau'+\epsilon,\cr
0 \qquad &\hbox{if}\quad t\geq \tau'+\epsilon.\enda\right.
\eeq
Define the Lipschitz function with compact support
 \bel{pen1}\varphi^{\epsilon}(t,y)~\doteq~\min\{ \rho^{\epsilon}(t,y),
\,\chi^\epsilon(t)\}.\eeq
Using $\vp^{\epsilon}$ as test function, since the family of measures
$\mu^t_-$ satisfy
the first equation in  (\ref{mb4}),
we obtain
\bel{vpe}
\int \bigg[ \int (\vp^{\epsilon}_t-c\vp^{\epsilon}_x) \, d\mu^t_- + \int
 \frac{c'}{2c}\,(R^2S - RS^2)\, \vp^{\epsilon}
\, dx\bigg] dt~=~0.
\eeq
We now observe that, if $\tau'$ is sufficiently close to $\tau$, then
for any $t\in [\tau, \tau']$ and $x$ close to $x^-(\tau)$, one has
$$0 ~=~ \vp^{\epsilon}_t + [-c(u(\tau, x(\tau)))+\ve_0 ] \vp^{\epsilon}_x ~
\leq~\vp^{\epsilon}_t -c( u(t,x)) \vp^{\epsilon}_x\,,$$
because  $-c(u(t,x))<-c(u(\tau, x(\tau)))+\ve_0 $ and $\vp^{\epsilon}_x\leq 0$.

Since the measures $\mu^t_-$ depend continuously on $t$
in the topology of weak convergence, taking the limit of (\ref{vpe}) as
 $\epsilon\to 0$, for $\tau,\tau'\notin\N$ we obtain
\bel{55}  \mu^{\tau'}_- \Big(]-\infty, \,X(\tau')]\Big)~
\geq~\mu^\tau_- \Big(]-\infty, \,X(\tau)]\Big)
+ \int_\tau^{\tau'}\int_{-\infty}^{X(t)}\frac{c'}{2c}[ R^2S-RS^2 ]\,dxdt\,.
\eeq
By (\ref{xpmt}), for  $t\in \,]\tau, \, \tau+\delta]$ one has
\bel{66}\mu_-^t\Big(\, ]-\infty,\, X(t)]\Big) ~\leq~\mu_-^t\Big(\, ]-\infty,\, x^-(t)[\Big) .\eeq
Using (\ref{55})-(\ref{66}) one obtains
\bel{ccc}\bega{l}\ds
\alpha(t) - \alpha(\tau) ~\geq~
\bigg[x^-(t) + \mu_-^t\Big( ]-\infty,\, x^-(t)[\Big)\bigg] - \bigg[ x^-(\tau) +
 \mu_-^\tau\Big( ]-\infty,\, x^-(\tau)]\Big)\bigg]\cr\cr
 \ds
\qquad \geq~ \Big[ - c\bigl(u(\tau, x^-(\tau))\bigr)+2\ve_0\Big](t-\tau) +
\int_\tau^t \int_{-\infty}^{x^-(s)}\frac{c'}{2c}[ R^2S-RS^2 ]\,dy ds + o(t-\tau).
\enda
\eeq
Since $\tau\notin \N$, the term
$$o(t-\tau) ~=~\int_\tau^t \int_{X(s)}^{x^-(s)}\frac{c'}{2c}[ R^2S-RS^2 ]\,dy ds
$$
is a higher order infinitesimal, namely
${o(t-\tau)\over t-\tau} \to 0$ as $t\to \tau+$.
Differentiating (\ref{ccc}) w.r.t.~$t$ at $t=\tau$, we obtain
$$\dot \alpha(\tau) ~\geq~ \Big[ - c\bigl(u(\tau, x^-(\tau))\bigr)+2\ve_0\Big]+
\int_{-\infty}^{x^-(\tau)}\frac{c'}{2c}[ R^2S-RS^2 ]\,dy ds$$
in contradiction
with (\ref{ceq5}).
\v
{\bf 5. } In this step we prove the uniqueness of the solution to
(\ref{ica})-(\ref{ceq5}).
Consider the weight
\bel{weight} W(t,\alpha)~\doteq~e^{\kappa A^+(t,\alpha)}
\eeq
with
\bel{A+}\qquad\qquad A^+(t,\alpha)~\doteq~\mu_+^t\Big(]-\infty,~x(t,\alpha)]\Big)+
[\zeta(T)-\zeta(t)].
\eeq
Here $\zeta$ is the function defined at (\ref{zdef}), while
\bel{mudef}\kappa~\doteq~\frac{M}{2c_0^2}\,.\eeq
 We recall that $\zeta(T)-\zeta(t)$
 provides an upper bound on the energy transferred from backward to forward
 moving waves and conversely, during the time interval $[t,T]$.
In turn, $A^+(t,\alpha)$ yields an upper bound on the total
energy of forward moving waves that can cross the backward characteristic
$x(\cdot,\alpha)$
during the time interval  $[t,T]$.

For any $\alpha_1<\alpha_2$ and $t\geq 0$, we  define a weighted distance
by setting
\bel{RD}
d^{(t)}(\alpha_1,\alpha_2)~\doteq~\int_{\alpha_1}^{\alpha_2}
W(t,\alpha)\, d\alpha.
\eeq

Consider two  solutions of (\ref{ceq5}), say
$\alpha_1(t)\leq\alpha_2(t)$.  For convenience, we use the shorter notation
$$x_i(t)~\doteq~x^-(t, \alpha_i(t)),\qquad\qquad i=1,2.$$
We recall that, by the definition of conservative solution,
the measures $c'(u)\cdot\mu^t_-$, and $c'(u)\cdot\mu^t_+$ are absolutely
continuous w.r.t.~Lebesgue measure for a.e.~time $t$.
As before, $R^2$ and $S^2$ denote the density of the absolutely continuous part of
$\mu^t_-$, and $\mu^t_+$, respectively. One has the estimate
\bel{dta2}
\bega{l}
\ds \int_{\alpha_1(t)}^{\alpha_2(t)}W(t,\alpha)\, d\alpha~-\int_{\alpha_1(\tau)}^{\alpha_2(\tau)}W(\tau,\alpha)\, d\alpha\cr\cr
 =~\ds \int_\tau^t\dot\alpha_2 W(s,\alpha_2) - \dot\alpha_1 W(s,\alpha_1)
+\int_{\alpha_1(t)}^{\alpha_2(t)}{\partial\over\partial t}W(s,\alpha)\, d\alpha ds\cr\cr
\ds =~\int_\tau^t\int_{\alpha_1(s)}^{\alpha_2(s)} \left\{{\partial\over \partial\alpha} \left[
\left(-c(u(s,x(s,\alpha)))+\int_{-\infty}^{x(s, \alpha)}
{c' (R^2S-RS^2)\over 2c}
\, dx\right) W(s,\alpha)\right]+{\partial\over\partial t}W(s,\alpha)
\right\}\, d\alpha ds\cr\cr
 =~\ds\int_\tau^t\int_{\alpha_1(s)}^{\alpha_2(s)}\left[ {c'(S-R)\over 2c( 1+R^2)}
+ {c' \,(R^2S-RS^2)\over 2c(1+R^2)}\right] W\, d\alpha
 \cr\cr
\qquad\qquad\ds+\int_{\alpha_1(s)}^{\alpha_2(s)} \left(-c
+\int_{-\infty}^{x(s, \alpha)}
{c' (R^2S-RS^2)\over 2c}
\, dx\right){S^2\over 1+R^2}\, \kappa W\, d\alpha ds
\cr\cr
\ds\qquad\qquad  - \int_\tau^t\int_{\alpha_1(s)}^{\alpha_2(s)}\left[ {1+2R^2\over 1+R^2} \,c S^2+
{1+R^2+S^2\over 1+R^2} \int_{-\infty}^{x(s,\alpha)}
{c' (R^2S-RS^2)\over 2c}\,dx + \dot \zeta(s)\right]\cdot\kappa W\, d\alpha ds\cr\cr
=~\ds\int_\tau^t\int_{\alpha_1(s)}^{\alpha_2(s)}\left[ {c'(S-R)\over 2c( 1+R^2)}
+ {c' \,(R^2S-RS^2)\over 2c(1+R^2)}\right] W\, d\alpha ds
-\int_\tau^t\int_{\alpha_1(s)}^{\alpha_2(s)} 2cS^2 \kappa W\, d\alpha ds
 \cr\cr
\qquad\qquad\ds  -\int_\tau^t\int_{\alpha_1(s)}^{\alpha_2(s)}  \left(  \int_{-\infty}^{x(s, \alpha)}
{c' (R^2S-RS^2)\over 2c}
\, dx+ \dot \zeta(s)\right)\, \kappa W\, d\alpha ds\cr\cr
\ds\leq ~\int_\tau^t\int_{\alpha_1(s)}^{\alpha_2(s)}\left\{\left\|c'\over 2c
\right\|_{\L^\infty} \Big(1+|S|+|S|+|S^2|\Big)- 2cS^2\kappa \right\} W\,d\alpha ds\cr\cr
\ds\leq ~\int_\tau^t\int_{\alpha_1(s)}^{\alpha_2(s)}\left\{
\left\|c'\over 2c\right\|_{\L^\infty}2(1+S^2)- 2cS^2\kappa \right\} W\,d\alpha ds~\leq~
\left\|c'\over 2c\right\|_{\L^\infty} \cdot \int_\tau^t\int_{\alpha_1(s)}^{\alpha_2(s)}
 W\,d\alpha ds\,.
 \enda
\eeq

By Gronwall's lemma this implies
\bel{Gro}
d^{(t)}\bigl(\alpha_1(t), \alpha_2(t)\bigr)~
\leq~ e^{C_0t}\, d^{(0)}\bigl(\alpha_1(0), \alpha_2(0)
\bigr),
\eeq
with $C_0\doteq \|c'/ 2c\|_{\L^\infty}$.
For every initial value $\bar \alpha$
this yields the uniqueness of the solution
of (\ref{ica})-(\ref{ceq5}).
\v
{\bf 6.} Finally, we claim that, for any initial data
(\ref{icy}), there exists a unique function $x^-(t)$
which satisfies the first equation in (\ref{dxpm}) together with (\ref{be1}).

Indeed, let
$x^-_1(t)$ and $x^-_2(t)$ be two  solutions
with $x^-_1(0)=x^-_2(0) = \bar y$.
For $i=1,2$, consider the functions
$$\alpha_i(t)~\doteq~x_i^-(t) + \int_{-\infty}^{\bar y}R^2(0,x)\, dx
+\int_0^t\int_{-\infty}^{x_i^-(t)}{c' \,(R^2S-RS^2)\over 2c}\, dx\, ds\,.$$
Then $x_i^-(t) = x(t,\alpha_i(t))$, and
both $\alpha_1,\alpha_2$ are solutions to the Cauchy problem (\ref{ceq4}),
with initial data
$$\alpha_1(0) ~=~ \alpha_2(0)~=~\bar y + \int_{-\infty}^{\bar y}R^2(0,x)\, dx\,.$$
The uniqueness result proved in step {\bf 5} now implies~
$x_1^-(t) = x(t,\alpha_1(t))=x(t,\alpha_2(t))=x_2^-(t)$.
\endproof
\v

\section{Lipschitz continuity in characteristic coordinates}
\setcounter{equation}{0}

Let $u=u(t,x)$ be a conservative solution to the wave equation (\ref{1.1})
with initial data (\ref{1.2}).
Given $(X,Y)\in\R^2$, there exists unique initial points
$\bar x= x_0(X)$ and $\bar y= y_0(Y)$ such such that
\bel{XY}
X~=~\bar x +\int_{-\infty}^{\bar x}R^2(0,x)dx,\qquad
\qquad Y~=~\bar y +\int_{-\infty}^{\bar y}S^2(0,x)dx.\eeq
By Lemma 2, there exists a unique backward characteristic $t\mapsto x^-(t,\bar x)$
starting at $\bar x$, and a unique forward characteristic $t\mapsto x^+(t,\bar y)$
starting at $\bar y$.  Indeed, recalling (\ref{xa})-(\ref{yb}),
we can write
\bel{xcar}x^-(t,\bar x) ~=~x(t, \alpha(t)),\qquad\qquad x^+(t,\bar y) ~=~y(t,\beta(t)),\eeq
where $\alpha(\cdot)$ provides a solution to (\ref{ceq5}) with initial data
$\alpha(0)=X$, and similarly for $\beta(\cdot)$.

Assuming that $\bar x\geq\bar y$, we define
$$P(X,Y)~=~(t(X,Y), x(X,Y))$$ to be the unique point
where these two characteristics cross.  That means
\bel{xtXY}x^-(t(X,Y),\,\bar x) ~=~x^+(t(X,Y),\,\bar y)~=~ x(X,Y).\eeq
We then define
\bel{uXY}
u(X,Y)~\doteq~u\bigl(t(X,Y), x(X,Y)\bigr).\eeq
\v
{\bf Lemma 3.} {\it The map $(X,Y)\mapsto (t,x,u)(X,Y)$ is  locally Lipschitz continuous.}
\v
{\bf Proof.}  {\bf 1.}
For a fixed $Y$, we show that the map
$X\mapsto (t,x,u)(X,Y)$ is locally Lipschitz continuous.
Fix an interval $[0,T]$ and let $X_1<X_2$.
Recalling the notation in (\ref{xa})-(\ref{yb}), consider the backward characteristics
\bel{bchar}t\mapsto x_1(t)~\doteq~(t, \alpha_1(t)),\qquad\qquad t\mapsto x_2(t)
~\doteq ~x(t, \alpha_2(t)),\eeq
where $\alpha_1,\alpha_2$ are the solutions of (\ref{ceq5}) with initial data
$\alpha_1(0)=X_1$ and $\alpha_2(0)=X_2$, respectively.
Similarly, let $t\mapsto y(t, \beta(t))$ be the forward characteristic, with $\beta(0)=Y$.

\begin{figure}[htbp]
\centering
\includegraphics[width=0.7\textwidth]{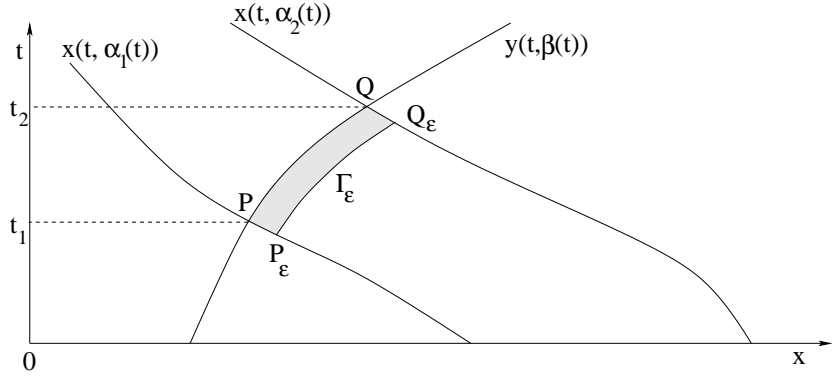}
\caption{Proving the Lipschitz continuity of $t,x,u$ as functions of
the variables $X,Y$.}
\label{f:wa38}
\end{figure}

As shown in Fig.~\ref{f:wa38}, assuming that $y(0, Y)\leq x(0, X_1)$, let
$t_1,t_2$ be the times when these characteristics  cross, so that
$x(t_i, \alpha_i(t_i))= y(t_i, \beta(t_i))$, $i=1,2$.
\v
{\bf 2.} According to (\ref{Gro}), we have
\bel{alip}\alpha_2(t_1)-\alpha_1(t_1)~\leq~C\, (\alpha_2(0)-\alpha_1(0))~=~C\,(X_2-X_1),\eeq
for some constant $C$ uniformly bounded as $t_1$ ranges over a bounded interval.
In turn this implies
\bel{xlip}x(X_2,Y)- x(X_1,Y)~\leq~x(t_1, \alpha_2(t_1))-x(t_1, \alpha_1(t_1))
~\leq~ \alpha_2(t_1)-\alpha_1(t_1)~\leq~C\,(X_2-X_1),\eeq
proving the Lipschitz continuity of the map $X\mapsto x(X,Y)$.

In turn, we have
\bel{tlip}t(X_2,Y)- t(X_1,Y)~\leq~\|c(u)\|_{\L^\infty}\cdot \bigl(x(X_2,Y)- x(X_1,Y)\bigr)
\leq~\|c(u)\|_{\L^\infty}\cdot C\,(X_2-X_1),\eeq
showing that the map $X\mapsto t(X,Y)$ is Lipschitz continuous as well.

{\bf 3.}
It remains to show that the map
$X\mapsto u(X,Y)$ is Lipschitz continuous.
Let $X_1<X_2$ and $Y$ be given.
As in step {\bf 1}, consider the backward characteristics
$t\mapsto x_i(t)\doteq x(t, \alpha_i(t))$, $i=1,2$, and the forward characteristic
$t\mapsto y(t)\doteq y(t, \beta(t))$, with  $\alpha_i(0)=X_i$ and
$\beta(0)=Y$.

As shown in Fig.~\ref{f:wa38}, consider the intersection points
$$P~=~(t_1, x(t_1, \alpha_1(t_1))), \qquad Q=(t_2, x(t_2, \alpha_2(t_2))). $$
Moreover, for $\eta>0$ small, consider the curve
$$t~\mapsto~\gamma_\eta(t)~\doteq~y(t, \beta(t)) + \eta\,,$$
and call $P_\eta$, $Q_\eta$ the points where this curve intersects
the two backward characteristics (\ref{bchar}).

Since $u=u(t,x)$ is continuous, one has
\bel{h}
u(X_2,Y) - u(X_1,Y) ~=~ u(Q)- u(P)~=~
\lim_{\ve\to 0+}~{1\over\ve}\int_0^\ve \Big(u(Q_\eta)- u(P_\eta)\Big)\, d\eta\,.\eeq
Call $\Gamma_\ve$ the rectangular region bounded by the two backward
characteristics in (\ref{bchar})  and by the curves $\gamma_0, \gamma_\ve$
(the shaded region in Fig.~\ref{f:wa38}).
We now compute
\bel{h1}\bega{rl}\ds\int_0^\ve u(Q_\eta)- u(P_\eta)\, d\eta&\ds=~
\dint_{\Gamma_\ve}
\Big[ u_t +c(u(t, y(t))) u_x - R(t,x)\Big]\, dx dt
+\dint_{\Gamma_\ve} R\, dx dt \cr\cr
&\ds\doteq~I_1+I_2\,.\enda
\eeq
We estimate the two integrals on the right hand side of (\ref{h1}).
Assuming $0\leq t_1<t_2\leq T$, by Cauchy's inequality we obtain
$$\bega{rl}|I_1|&\ds\leq~\dint_{\Gamma_\ve}\Big|c(u(t,y(t)) - c(u(t,x))\Big|
\,|u_x(t,x)|\, dxdt
\cr\cr
&=~\ds \O(1)\cdot \int_0^T \left(\int_{0<x-y(t)<\ve}|x-y(t)|^{1/2}\, |u_x(t,x)|\, dx
\right)dt
\cr\cr
&=~\ds \O(1)\cdot\int_0^T \ve \left(\int_{0<x-y(t)<\ve}|u_x(t,x)|^2\, dx
\right)^{1/2}dt\,.
\enda
$$

This implies
\bel{h2}
\lim_{\ve\to 0+} {1\over \ve} ~I_1~=~0.\eeq

Next,
\bel{h3} |I_2|~\leq~\int_{\Gamma_\ve} (1+R^2)\, dx dt\eeq
To estimate the above  integral, consider the weight function
$$V^\ve(t,x)~\doteq~\left\{ \bega{cl} 0\quad & \hbox{if} \quad x \leq y(t),\cr\cr
\ds{x-y(t)\over\ve}~\quad & \hbox{if} \quad x \in
 [ y(t),~ y(t) +\ve],
\cr\cr
1\quad & \hbox{if} \quad x \geq y(t) +\ve\,.
\enda\right.
$$
Since the family of measures $\mu^t_-$ satisfies the balance equation in (\ref{mb4})
with speed $<-c_0$, for a.e.~time $t$ such that
$x_1(t)~\le ~y(t)~\le~ y(t)+\ve ~\le ~x_2(t)$ we have
$$\bega{rl}\ds{1\over\ve}\,\int_{y(t)}^{y(t) +\ve} R^2(t,x)\, dx
&\ds\leq~{1\over\ve}\, \mu^t_-\bigl(\,]x_1(t),~x_2(t)[\,\bigr)\cr\cr
&\ds
\leq~-{1\over 2c_0}\,{d\over dt}
\left[\int_{x_1(t)}^{x_2(t)} V(t,x) d\mu^t_-\right]
+\int_{x_1(t)}^{x_2(t)}{|c'(u)|\,\bigl|R^2 S - R S^2\bigr|\over 2c}\, dx.
\enda $$
Integrating over the interval $[t_1,t_2]$ and using the estimate (\ref{dta2})
to control the total contribution of  the source term,    we obtain
\bel{h4}{|I_2|\over\ve}~\leq~{1\over\ve}\int_{\Gamma_\ve} (1+R^2)\, dx dt
~\leq~(t_2-t_1) + {1\over 2c_0}(X_2-X_1) + \O(1)\cdot (X_2-X_1).\eeq
This proves the Lipschitz continuity of the map $X\mapsto u(X,Y)$.

  The Lipschitz continuity of $(t,x,u)$ as functions of $Y$ is proved in exactly the same way.
\endproof
\v
{\bf Remark 3.} By Rademacher's theorem, the above result implies that
the map
\bel{Ldef}\Lambda:(X,Y)\mapsto (t(X,Y), \, x(X,Y))\eeq
 is a.e.~differentiable.   We can thus consider the set
$\Omega$ of {\it critical points}  and the set $V$
 of {\it critical values} of $\Lambda$, by setting
\bel{critp}
\Omega~\doteq~\Big\{ (X,Y)\,;~~\hbox{either $D\Lambda(X,Y)$ does not exists,
or else det$D\Lambda(X,Y)=0$}\Big\}.\eeq
\bel{critv}
V~\doteq~\Big\{ \Lambda(X,Y)\,;~~(X,Y)\in \Omega\Big\}.\eeq
By the area formula \cite{Zi},  the 2-dimensional measure of $V$ is zero.

In general, the map $\Lambda:\R^2\mapsto\R^2$ is onto but not
one-to-one.   However,
for each $(t_0,x_0)\notin V$, there exist a unique point $(X,Y)$ such that
$\Lambda(X,Y) = (t_0, x_0)$.

Next, consider a  function $f(t,x)$, with  $f\in\L^1(\R^2)$.
Since $f$ is defined up to a set of measure zero in the $t$-$x$ plane,
the composition
$\tilde f(X,Y) = f(\Lambda(X,Y))$
is well defined at a.e.~point $(X,Y)\in \R^2\setminus\Omega$.
Moreover, we have the change of variable formula
\bel{icv}\int_{\R^2} f(t,x)\, dxdt~=~\int_{\R^2\setminus\Omega} \tilde f(X,Y)
\cdot |\det D\Lambda(X,Y)|\, dXdY.\eeq
To compute the determinant of the Jacobian matrix $D\Lambda$, we observe
that
\bel{pdXY}x_X~=~c(u) \, t_X\,,\qquad\qquad x_Y~=~-c(u) \, t_Y\,,\eeq
\bel{dL}
 D\Lambda~=~\left(\bega{ccc} t_X & & t_Y\cr
\cr
x_X&& x_Y  \enda\right) ~=~
\left(\bega{ccc} {x_X\over c(u)}&& -{x_Y\over c(u)}   \cr
\cr
x_X&& x_Y\enda\right).\eeq
Hence
\bel{detL} \bigl|\det D\Lambda\bigr|~=~ {2\over c(u)}
x_X x_Y\,.\eeq


For future use, in the $X$-$Y$ plane we define the ``good set"
\bel{gset}\G ~\doteq~\R^2\setminus\Omega\,.\eeq

\v
\section{An equivalent semilinear system}
\setcounter{equation}{0}

In this section we introduce further variables and show that, as functions of
$X,Y$, these variables satisfy a semilinear system with smooth coefficients.
In particular, their values are uniquely determined by the initial data.
By showing that the map $(X,Y)\mapsto (t,x,u)(X,Y)$ is uniquely determined,
we eventually obtain the uniqueness of the solution $u(t,x)$ of the Cauchy problem
(\ref{1.1})-(\ref{1.2}).
\v
 Recalling (\ref{xadef})-(\ref{ybdef}),
for given  initial values $\bar \alpha, \bar \beta$
let  $t\mapsto \alpha(t,\bar \alpha)$ and $t\mapsto \beta(t,\bar \beta)$ be the unique
solutions to
$$\bega{rl} \alpha(t)&\ds=~\bar\alpha + \int_0^t \left(-c(0)+\int_{-\infty}^{x(t, \alpha(t))}
{c' (S-R+R^2S-RS^2)\over 2c}
\, dx\right)dt\,,\cr\cr
\beta(t)&\ds=~\bar\beta + \int_0^t \left(c(0)-\int_{-\infty}^{y(t, \beta(t))}
{c' (S-R+R^2S-RS^2)\over 2c}
\, dx\right)dt\,.\enda$$
The existence and uniqueness of these functions was proved in Section~3.
We recall that $t\mapsto x^-(t) = x(t,\alpha(t))$ and $t\mapsto x^+(t)
=y(t,\beta(t))$ are then the  unique  backward  and forward characteristics
starting from the points $x(0, \bar\alpha)$ and $y(0, \bar\beta)$, respectively.
Define the new dependent variables $p(X,Y)$ and $q(X,Y)$ by setting
\bel{pdef}
p(X,Y) ~=~{\partial \over \partial \bar\alpha} \alpha(\tau,\bar\alpha)
\bigg|_{\bar\alpha=X,\, \tau = t(X,Y)}\,,\qquad
q(X,Y)~=~{\partial \over \partial \bar\beta} \beta(\tau,\bar\beta)\bigg|_{\bar\beta=Y,
\, \tau = t(X,Y)}\,.\eeq
In addition, recalling the definitions of the maps $\alpha\mapsto x(t,\alpha)$
and $\beta\mapsto y(t,\beta)$ in (\ref{xadef})-(\ref{ybdef}), we define
\bel{nedef}
\nu(X,Y)~\doteq~\frac{\partial x}{\partial \alpha}(t(X,Y),\alpha(t,x(X,Y))),\qquad\qquad \eta(X,Y)~\doteq~\frac{\partial x}{\partial \beta}(t(X,Y),\beta(t,x(X,Y)))\,.\eeq
Finally, observing that the functions $c, p,q$ are strictly positive, we define
$\xi,\zeta$ by setting
\bel{tvt}
\xi(X,Y) ~\doteq~{2c(u(X,Y))\over p(X,Y)} u_X(X,Y)\,,\qquad\qquad
\zeta(X,Y) ~\doteq~{2c(u(X,Y))\over q(X,Y)} u_Y(X,Y)\,.\eeq
By Rademacher's theorem, the above derivatives are a.e.~well defined, because
\begi
\item[(i)] for any $t$, the functions $\bar\alpha\mapsto \alpha(t,\bar \alpha)$,
$\bar\beta\mapsto \beta(t,\bar \beta)$, $\alpha\mapsto x(t,\alpha)$,  and $\beta\mapsto y(t,\beta)$ are Lipschitz continuous, and
\item[(ii)] by Lemma~3, both $x$ and $u$ are Lipschitz continuous functions of $X,Y$.
\endi
Moreover,
$$p(X,Y) ~=~q(X,Y)~=~1\qquad\hbox{if}\qquad t(X,Y)~=~0.$$
Our main goal is to show that these variables satisfy the semilinear system
with smooth coefficients
\bel{sls}
\left\{\bega{l}
u_X=\frac1{2c}\xi p,\qquad
u_Y=\frac1{2c}\zeta q,\cr\cr
x_X=\frac12 \nu p,\quad  \quad x_Y=-\frac12\eta q,\cr\cr
t_X=\frac{1}{2c}\nu p,\quad  \quad  t_Y=\frac1{2c}\eta q,\cr\cr
p_Y=\frac{c'}{4c^2}[\zeta-\xi]pq,\cr\cr
q_X=\frac{c'}{4c^2}[\xi-\zeta]pq,\cr\cr
\nu_Y=\frac{c'}{4c^2}\xi(\nu-\eta)q,\cr\cr
\eta_X=-\frac{c'}{4c^2}\zeta(\nu-\eta)p,\cr\cr
\xi_Y=-\frac{c'}{8c^2}(\eta+\nu)q+\frac{c'}{4c^2}(\xi^2+\eta\nu)q\,,\cr\cr
\zeta_X=-\frac{c'}{8c^2}(\eta+\nu)p+\frac{c'}{4c^2}(\zeta^2+\eta\nu)p\,.
\enda\right.
\eeq
More precisely, we have
\v
{\bf Theorem 4.} {\it By possibly changing the functions $p,q,\nu,\eta,\xi,\zeta$
on a set of measure zero in the $X$-$Y$ plane, the following holds.
\begi
\item[(i)] For a.e.~$X_0\in \R$, the functions $t,x,u,p,\nu,\xi$ are absolutely continuous on every vertical segment of the form $S_0\doteq \{ (X_0,Y)\,;~~~a<Y<b\}$.   Their partial derivatives w.r.t.~$Y$
satisfy a.e.~the corresponding equations in (\ref{sls}).

\item[(ii)] For a.e.~$Y_0\in \R$, the functions $t,x,u,q,\eta,\zeta$ are
absolutely continuous on every horizontal segment of the form
$S_0\doteq \{ (X,Y_0)\,;~~~a<X<b\}$.   Their partial derivatives w.r.t.~$X$
satisfy a.e.~the corresponding equations in (\ref{sls}).
\endi}
\v
Toward a proof, we recall a standard result in the theory of Sobolev spaces.
\v
{\bf Lemma 4.}  {\it  Let $\Gamma = \, ]a,b[\, \times\, ]c,d[\,$
be a rectangle  in the $X$-$Y$ plane.  Assume that $u\in \L^\infty (\Gamma)$
has a weak partial derivative w.r.t.~$X$. That means
\bel{fu}
\int_\Gamma (u\vp_X + f\vp)\, dXdY~=~0\eeq
for some $f\in \L^1(\Gamma)$ and all test functions $\vp\in \C^\infty_c(\Gamma)$.
Then, by possibly modifying $u$ on a set of measure zero, the following holds.
For a.e.~$Y_0\in \,]c,d[\,$, the map $X\mapsto u(X, Y_0)$ is absolutely continuous
and
\bel{pdX}{\partial \over\partial X} u(X,Y_0)~=~f(X, Y_0)\qquad\qquad\hbox{for a.e.~} X\in \,]a,b[\,.
\eeq}
\v
For a proof, see for example \cite{BFA}, p.159, or \cite{Zi}, p.44.
To use the above result, it is convenient to replace the test functions $\vp$
with characteristic functions of arbitrary rectangles contained in $\Gamma$.
\v
{\bf Lemma 5.} {\it   Let $\Gamma = \, ]a,b[\, \times\, ]c,d[\,$
be a rectangle  in the $X$-$Y$ plane.  Assume that $u\in \L^\infty (\Gamma)$
and $f\in \L^1(\Gamma)$.
Moreover assume that there exists null sets $\N_X\subset\,]a,b[\,$ and
$\N_Y\subset\,]c,d[\,$ such that the following holds.

For every $\ov X_1,\ov X_2\notin \N_X$ and $\ov Y_1,\ov Y_2\notin \N_Y$ with $\ov X_1<\ov X_2$ and
$\ov Y_1<\ov Y_2$, one has
\bel{IXY}
\int_{\ov Y_1}^{\ov Y_2} \Big[ u(\ov X_2,Y) - u(\ov X_1,Y)\Big]\, dY ~=~\int_{\ov Y_1}^{\ov Y_2}
\int_{\ov X_1}^{\ov X_2} f(X,Y)\, dXdY\,.
\eeq
Then the conclusion of Lemma 4 holds.}
\v
{\bf Proof.} Consider any test function $\vp\in \C^\infty_c(\Gamma)$.
We need to show that (\ref{fu}) holds.
Given $\ve>0$, we can find points
$$a=X_0<X_1<\cdots<X_N=b\,,\qquad\qquad c=Y_0<Y_1<\cdots<Y_N=d$$
with $X_i\notin\N_X$, $Y_i\notin \N_Y$ and
$$X_i-X_{i-1}~<~\ve\,,\qquad Y_i-Y_{i-1}~<~\ve\,,\qquad\quad  i=1,\ldots,N\,.$$
Define the approximate function $\vp^\ve$ by setting
$$\vp^\ve(X,Y)~\doteq~\vp(X_{i-1}, Y_{i-1}) + (X-X_{i-1}) \cdot \bigl[ \vp(X_i, Y_{i-1})-
\vp(X_{i-1}, Y_{i-1}) \bigr].$$
Taking a sequence of these approximations with  $\ve\to 0$, we have the convergence
$$\|\vp^\ve-\vp\|_{\L^\infty} ~\to~ 0,\qquad\qquad \|\vp^\ve_X-\vp_X\|_{\L^\infty} ~\to~ 0.
$$
Therefore,
\bel{fu2}
\int_\Gamma (u\vp_X + f\vp)\, dXdY~=~\lim_{\ve\to 0}
\int_\Gamma (u\vp^\ve_X + f\vp^\ve)\, dXdY ~=~0.\eeq
Indeed, for every approximate function $\vp^\ve$ one has
$$\bega{l}
\ds \int_\Gamma (u\vp^\ve_X + f\vp^\ve)\, dXdY ~=~\sum_{i,j=1}^N
\int_{X_{i-1}}^{X_i}\int_{Y_{j-1}}^{Y_j} (u\vp^\ve_X + f\vp^\ve)\, dXdY
\cr\cr
\ds \quad =~\sum_{i,j=1}^N \int_{Y_{j-1}}^{Y_j}
\int_{X_{i-1}}^{X_i}\left[ u(X_{i-1},Y) + \int_{X_{i-1}}^X f(X',Y)\, dX'\right]\,\vp^\ve_X\,
dXdY\cr\cr
\qquad\qquad \ds  +\sum_{i,j=1}^N
\int_{Y_{j-1}}^{Y_j}\int_{X_{i-1}}^{X_i} f\vp^\ve\, dXdY \cr\cr
\ds \quad =~- \sum_{i,j=1}^N
\int_{Y_{j-1}}^{Y_j}\int_{X_{i-1}}^{X_i} f\vp^\ve\, dXdY \cr\cr
\qquad\qquad \ds
+ \sum_{i,j=1}^N
\int_{Y_{j-1}}^{Y_j} \Big[ u(X_i,Y) \vp(X_i, Y_{i-1}) - u(X_{i-1},Y)\vp(X_{i-1} Y_{i-1})\Big]\, dY\cr\cr
\qquad\qquad \ds  +\sum_{i,j=1}^N
\int_{Y_{j-1}}^{Y_j}\int_{X_{i-1}}^{X_i} f\vp^\ve\, dXdY\cr\cr \quad =~0\,.
\enda$$
\endproof

\vsk
Recalling the sets of regular and critical points
$\G $, $\Omega$,  defined at (\ref{gset}) and (\ref{critp}) respectively, we now
derive a representation for the  variables $\nu,\eta,\xi,\zeta$ in terms of $R$ and $S$.

{\bf Lemma 6.} {\it
\begi
\item[(i)] If $(X,Y)\in \G ~\doteq~\R^2\setminus\Omega$ then
\bel{L40} {p(X,Y)\over x_X(X,Y)} ~=~2(1+R^2)\,,\qquad\qquad
{q(X,Y)\over x_Y(X,Y)}~=~2(1+S^2)\,,\eeq
\bel{L41}
\left\{\bega{rl}\nu(X,Y)&=~\frac1{1+R^2}\,,\cr\cr  \eta(X,Y)&=~\frac1{1+S^2}\,,
\enda\right.\qquad\qquad
\left\{\bega{rl}\xi(X,Y)&=~\frac{R}{1+R^2}\,,\cr\cr  \zeta(X,Y)&=~\frac{S}{1+S^2}\,,
\enda\right.
 \eeq
where the right hand sides are evaluated at the point $(t(X,Y),x(X,Y))$.

\item[(ii)] For a.e.~$(X,Y)\in\Omega$, one has
\bel{L42}\nu(X,Y)~=~ \eta(X,Y)~=~\xi(X,Y)~=~\zeta(X,Y)~=~0\,.\eeq

\endi
}
\v
{\bf Proof.}
{\bf 1.} Consider a regular value $(X,Y)\in \G $.  To fix the ideas,
let $t(X,Y)=\tau$. Recalling the definition (\ref{xadef}) and the fact that
the absolutely continuous part of $\mu_-^\tau$ has density $R^2$, we conclude
\bel{90}{\partial  \over\partial \alpha}x(\tau,\alpha) ~=~\frac1{1+R^2}.\eeq
On the other hand,
\bel{91}{\partial \over\partial \bar\alpha} x(\tau, \alpha(\tau,\bar\alpha))
~=~{2x_X}\,.\eeq
Together, the above equalities yield the first identity in (\ref{L40}). The
second one is proved similarly.
\v
{\bf 2.} The first identity in
(\ref{L41}) is precisely (\ref{90}), and the second one is similar.
To prove the third identity we observe that, at a point $(X,Y)\in \G$,
$$u_X~=~(u_t + c(u) u_x)\, t_X~=~R\, {x_X\over c(u)}\,.$$
Therefore, by (\ref{L40}),
$$\xi~\doteq ~{2c(u)\over p}\cdot u_X~=~{2c(u)\over p}\cdot R\, {x_X\over c(u)}~=~
{R\over 1+R^2}\,.$$
The last identity in (\ref{L41}) is proved similarly.
\v
{\bf 3.} Finally, consider the set $\G_0$ of points $(X,Y)$ where the map $\Lambda:(X,Y)\mapsto
(t,x)$ is differentiable but either $x_X=0$ or $x_Y=0$.
By Rademacher's theorem, the set $\R^2\setminus (\G\cup\G_0)$ has measure zero.

Assume $(X,Y)\in \G_0$, with $x_X(X,Y)=0$.
Then
$$\nu(X,Y)~=~{\partial x\over\partial \bar\alpha}\cdot
\left({\partial\alpha\over\partial \bar\alpha}\right)^{-1}~=~0.$$
In the same way one proves that  $\eta=0$.

Next, we claim that $t_X(X,Y)=x_X(X,Y)=0$ implies $u_X=0$.
This will be proved by refining the
estimates in step~{\bf 3} of the proof of Lemma~3.   Adopting the same construction,
for any $\delta>0$ we can replace the bound (\ref{h3}) with
\bel{h31}|I_2|~\leq~\int_{\Gamma_\ve}(C_\delta + \delta R^2) \, dxdt\,,\eeq
where $C_\delta = (4\delta)^{-1}$.
In this way, the estimate (\ref{h4}) can be replaced by
\bel{h41}{|I_2|\over\ve}~\leq~{1\over\ve}\int_{\Gamma_\ve} (C_\delta+\delta R^2)\, dx dt
~\leq~C_\delta (t_2-t_1) + \delta\left[{1 \over 2c_0}(X_2-X_1) + \O(1)\cdot (X_2-X_1)
\right].\eeq
Letting $\ve\to 0$ this implies
$$|u(X_2,Y)- u(X_1,Y)|~\leq~C_\delta \bigl|t(X_2,Y)-t(X_1,Y)\bigr| + \delta\left[{1 \over 2c_0}(X_2-X_1) + \O(1)\cdot (X_2-X_1)
\right]$$
In the above formula we can now take $X_1=X$, $X_2= X+\ve$.   Letting $\ve\to 0$ one obtains
$$\bigl|u_X(X,Y)\bigr|~=~\left|\lim_{\ve\to 0} {u(X+\ve,Y)-u(X,Y)\over\ve}\right| ~
\leq~C_\delta t_X + \delta \cdot \O(1)\,.$$
Since $t_X(X,Y)=0$ and $\delta>0$ can be taken arbitrarily small, this proves $u_X(X,Y)=0$.
From the definition (\ref{tvt}) it thus follows $\xi(X,Y)=0$.
Similarly, if $t_Y=x_Y=0$, then $\zeta =0$.
\endproof
{\bf Remark 4.} By a direct calculation, one finds
 \bel{detDL} \det D\Lambda~=~\frac{pq}{2c(1+R^2)(1+S^2)}.
\eeq
\v

\subsection{Proof of Theorem 4.}
To achieve a proof of Theorem 4, we will show that the assumptions of Lemma~5
apply to all the  variables $t,x,u,p,q,\eta,\nu,\xi,\zeta$ in \eqref{sls}.
For thus purpose, consider any rectangle $\Q\doteq
[X_1,X_2]\times [Y_1,Y_2]$ in the $X$-$Y$ plane.
\v
{\bf 1 - Equations for $u$.}  By Lemma~3, the function $u$ is Lipschitz continuous
w.r.t.~the variables $X,Y$.   The equations
$$u_X~=~{1\over 2c} \, \xi p\,,\qquad\qquad u_Y~=~{1\over 2c} \, \zeta q$$
are immediate consequences of the definitions (\ref{tvt}).
\v
{\bf 2 - Equations for $x$ and $t$.}   By Lemma~3, both functions $x,t$
are Lipschitz continuous
w.r.t.~the variables $X,Y$. Recalling the definitions (\ref{nedef})-(\ref{tvt}), we compute
\beq\bega{ll}
\ds\frac{\partial}{\partial X}x(X,Y)&\ds=\frac12\frac{\partial}{\partial\bar\alpha}x(t(X,Y),\alpha(t(X,Y),\bar\alpha))\cr\cr
&\ds=\frac12\frac{\partial x}{\partial \alpha}(t(X,Y),\alpha(t(X,Y),x(X,Y)))\frac{\partial\alpha}{\partial\bar\alpha}(t(X,Y),\bar\alpha)\bigg|_{\bar\alpha=X}\cr\cr
&\ds=\frac12\nu p.
\enda
\eeq
The equation for $x_Y$ is obtained in a similar way.

In turn, the equations for $t$ are derived by
\beq
t_X~=~\frac{x_X}{c(u)}~=~\frac{1}{2c}\nu p,\qquad \qquad t_Y~=~-\frac{x_Y}{c(u)}~=~
\frac1{2c}\eta q.
\eeq

\v
{\bf 3 - Equations for $p$ and $q$.}
By Lemma~3, the map $\Lambda(X,Y) \doteq  (t(X,Y),\, x(X,Y))$
is Lipschitz continuous, hence its Jacobian matrix
$D\Lambda$ is a.e.~well defined.
Consider the domain
\bel{Dedef}
\D~\doteq~\Big\{ (X,Y)\,;\quad X\in  [X_1, \, X_2],\quad  Y\in [Y_1,Y_2],\qquad
 \hbox{det}~D\Lambda(X,Y) >0\Big\}.\eeq
Recalling  (\ref{ceq4}) and using (\ref{icv}), we obtain
\beq
\bega{l}
\ds\int_{X_1}^{X_2}p(X,Y_2)-p(X,Y_1)dX\cr\cr
\qquad =~\ds\int_{X_1}^{X_2} \left[\frac{\partial\alpha(\tau,X)}{\partial X}
\bigg|_{\tau=\tau(X,Y_2)}-
\frac{\partial\alpha(\tau,X)}{\partial X}\bigg|_{\tau=\tau(X,Y_1)}\right]dX\cr\cr
\qquad =~\ds\int_{X_1}^{X_2} \left[\frac{\partial}{\partial X}\int_{x(X,Y_1)}^{x(X,Y_2)}\int_{\tau(\Tilde  X,Y_1)}^{\tau(\Tilde X,Y_2)} \frac{c'}{2c}(S-R+R^2S-S^2R)dtdx\right]d\Tilde X\cr\cr
\qquad =~\ds\iint_{\Lambda(\D)}\frac{c'}{2c}(S-R+R^2S-S^2R)dxdt\cr\cr
\qquad\ds=~\iint_{\D}
\frac{c'}{2c}(S-R+R^2S-S^2R)\cdot \hbox{det} D\Lambda(X,Y)\, dXdY\cr\cr
 \qquad\ds=~\iint_{\D}
\frac{c'}{4c^2}(S-R+R^2S-RS^2)\frac1{1+R^2}\frac1{1+S^2}pq \,dXdY\cr\cr
\qquad=~\ds\iint_\Q\frac{c'}{4c^2}(\zeta-\xi)pq\,dXdY.
\enda
\eeq
The last equality follows from Lemma~6, part (i) for the integral over $\D$ and
part (ii) for the integral over $\Q\setminus \D$.

Thus by Lemma 4 and 5,
\bel{pYeq}
p_Y~=~\frac{c'}{4c^2}(\zeta-\xi)pq.
\eeq
Similarly,
$$
q_X~=~\ds-\frac{c'}{4c^2}(\zeta-\xi)pq.
$$
\v
{\bf 4 - Equations for $\eta$ and $\nu$.}
We first observe that
$$
\bega{l}
\ds\int_{X_1}^{X_2}\Big[p\nu(X,Y_2)-p\nu(X,Y_1)\Big]\,dX
~ =~\ds\int_{X_1}^{X_2} \left[ \frac{\partial x(\tau,X)}{\partial X}
\bigg|_{\tau=t(X,Y_2)}-
\frac{\partial x (\tau,X)}{\partial X}\bigg|_{\tau=t(X,Y_1)}\right]\,dX~
\cr\cr
\qquad =~\ds\int_{X_1}^{X_2} \left[\frac{\partial}{\partial X}
\int_{x(X,Y_1)}^{x(X,Y_2)}\int_{t(\Tilde X,Y_1)}^{t(\Tilde X,Y_2)} \frac{c'}{2c}(S-R)dtdx\right]\,d\Tilde X
~ =~\ds\iint_{\Lambda(\D)}\frac{c'}{2c}(S-R)\,dxdt\,.\cr\cr
\enda
$$
In turn, by Remark 3 we obtain
\bel{eqp}
\bega{ll}
\ds\int_{X_1}^{X_2}\Big[ p\nu(X,Y_1)-p\nu(X,Y_2)\Big]\,dX&\ds=~\int_{\D}
\frac{c'}{2c}(S-R)\cdot \hbox{det} D\Lambda(X,Y)\, dXdY\cr\cr
  &\ds=~\int_{\D}
\frac{c'}{4c^2}(S-R)\frac1{1+R^2}\frac1{1+S^2}pq \,dXdY\cr\cr
&=~\ds\int_{\Q}\frac{c'}{4c^2}(\nu\zeta-\xi\eta)pq\, dXdY
\enda
\eeq
By Lemma 4 and 5, the above implies
\beq
(p\nu)_Y~=~\frac{c'}{4c^2}(\nu\zeta-\xi\eta)pq.
\eeq
Recalling the equation (\ref{pYeq}) for $p$,
we obtain
\beq
\nu_Y=\frac{c'}{4c^2}\xi(\nu-\eta)q.
\eeq
Similarly,
\beq
\eta_X=-\frac{c'}{4c^2}\zeta(\nu-\eta)p.
\eeq
\v
\begin{figure}[htbp]
\centering
\includegraphics[width=1\textwidth]{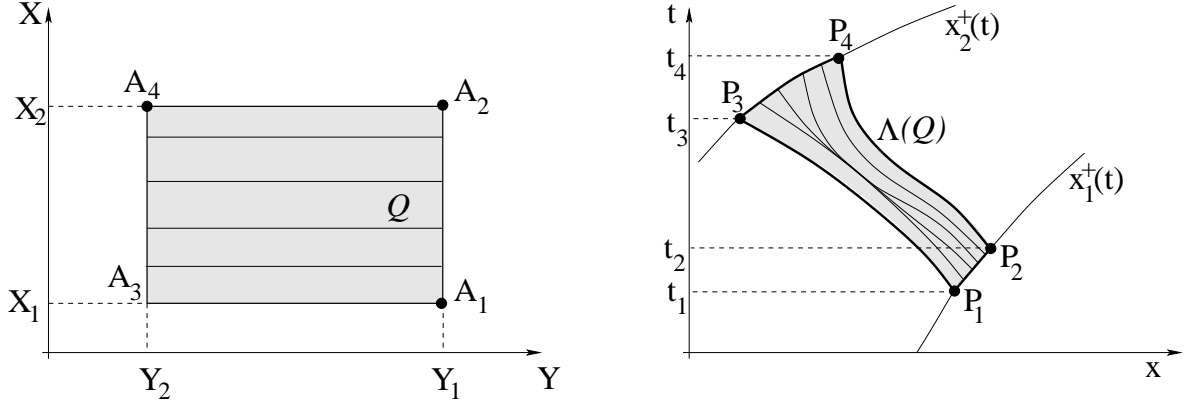}
\caption{\small Left: the rectangle $\Q$ in the $X$-$Y$ plane.
Right: the image $\Lambda(\Q)$ in the $t$-$x$ plane. For
$i=1,2,3,4$, the points $P_i=\Lambda(A_i)$ are defined as in
 (\ref{Pidef}).}
\label{f:hyp201}
\end{figure}
{\bf 5 - Equations for $\xi$ and $\zeta$.}
We observe that, by (\ref{2.2})-(\ref{2.3}),
$R$ provides a weak solution to the balance law
\bel{blR}R_t - (cR)_x~=~{c'\over 4c}(R^2-S^2) - c' u_x R~=~
{c'\over 4c}(R^2-S^2) - {c'\over 2c}  (R-S)R=-\frac{c'}{4c}(R-S)^2.
\eeq
Notice that, by definition of conservative solution,
the right hand side is a function in $\L^1(\R^2)$, w.r.t.~the variables $t$, $x$.

Next, we wish to characterize the distributional derivative $u_{XY}$.
More precisely, we seek a function $f\in \L^1_{loc}(\R^2)$ such that the following holds.
Consider any values $X_1<X_2$ and $Y_1> Y_2$.
Then
\bel{fdef}
[u(X_2, Y_1) - u(X_1, Y_1)]-[u(X_2, Y_2) - u(X_1, Y_2)] ~
=~\int_{X_1}^{X_2}\int_{Y_2}^{Y_1} f(X,Y)\, dX dY.\eeq
Toward this goal,
  consider the image of these four points under the map $\Lambda$, in the $t$-$x$ plane:
\bel{Pidef}\bega{l}
P_1\doteq (t_1,x_1)=\Lambda(X_1,Y_1), \qquad P_2\doteq(t_2,x_2)=\Lambda(X_2,Y_1),\cr\cr
P_3\doteq (t_3,x_3)=\Lambda(X_1,Y_2), \qquad P_4\doteq(t_4,x_4)=\Lambda(X_2,Y_2),
\enda\eeq
  (see Fig.~\ref{f:hyp201}).
We now construct a family of test functions $\phi^\epsilon$ approaching the
characteristic function of the set $\Lambda(\Q)$, where $\Q\doteq [X_1,X_2]\times [Y_1,Y_2]$.
More precisely:
 \bel{pen11}\phi^{\epsilon}(s,y)~\doteq~\min\{ \varrho^{\epsilon}(s,y),
\,\varsigma^\epsilon(s,y)\},\eeq
where
\bel{testf1}\varrho^{\epsilon}(s,y)~\doteq~\left\{\bega{cl} 0 \qquad &\hbox{if}\quad y~
\leq ~ x_1^-(s)-\epsilon\cr
1+\epsilon^{-1}(y-x_1^-(s)) \qquad &\hbox{if}\quad x_1^-(s)-\epsilon~\leq~ y\leq
x_1^-(s)\cr
1\qquad &\hbox{if}\quad x_1^-(s)~ \leq ~y~\leq~ x_2^-(s)\cr
1-\epsilon^{-1}(y-x_2^-(s))\qquad &\hbox{if}\quad  x_2^-(s)~\leq ~
y~\leq~ x_2^-(s)+\epsilon\cr
0\qquad &\hbox{if}\quad y~\geq ~x_2^-(s)+\epsilon\,,\enda\right.\eeq
\bel{testf2}\varsigma^\epsilon(s,y)~\doteq~\left\{\bega{cl} 0 \qquad &\hbox{if}\quad y
~\leq~  x_1^+(s)-\epsilon\cr
1+\epsilon^{-1}(y-x_1^+(s)) \qquad &\hbox{if}\quad x_1^+(s)-\epsilon~\leq~ y~\leq~
x_1^+(s)\cr
1\qquad &\hbox{if}\quad x_1^+(s)~ \leq~ y~\leq~ x_2^+(s)\cr
1-\epsilon^{-1}(y-x_2^+(s))\qquad &\hbox{if}\quad  x_2^+(s)~\leq~ y~\leq ~
x_2^+(s)+\epsilon
\cr
0\qquad &\hbox{if}\quad y~\geq~ x_2^+(s)+\epsilon\,.\enda\right.\eeq
Here $t\mapsto x_1^-(t)$ and $t\mapsto x_2^-(t)$ are the backward characteristics
corresponding to $X=X_1$ and $X=X_2$ respectively.   Similarly,
$t\mapsto x_1^+(t)$ and $t\mapsto x_2^+(t)$ are the forward characteristics
corresponding to $Y=Y_1$ and $Y=Y_2$ respectively.

\begin{figure}[htbp]
\centering
\includegraphics[width=1\textwidth]{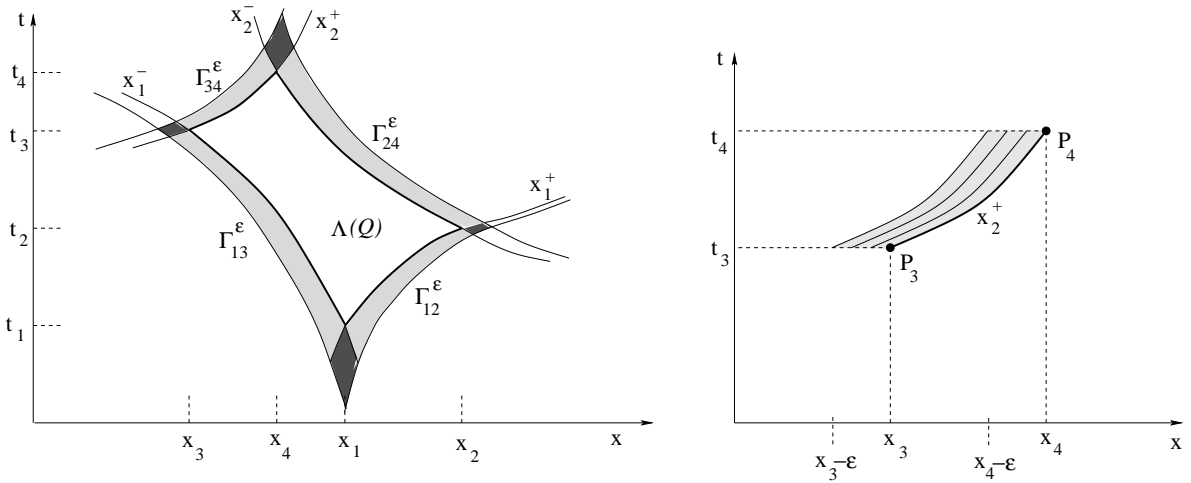}
\caption{\small   Left: the support of the test function $\phi^\epsilon$ in (\ref{pen11}).
Right: the difference $u(P_4)-u(P_3)$ can be approximated
by the integral of a directional derivative of $u$
over the shaded area.}
\label{f:hyp203}
\end{figure}

By (\ref{2.3}) for every test function $\varphi\in \C^1_c(\R^2)$ we now have
$$
\iint R\Big[\varphi_t -(c\varphi)_x\Big] \,dxdt~=~
-\iint \frac{c'}{4c}(R^2-S^2)\varphi\, dxdt\,.
$$
Since
$$R\in \L^2_{loc}(\R^2),\qquad  \quad c(u)\in H^1_{loc}(\R^2),
\qquad\quad {1\over c(u)}\in  H^1_{loc}(\R^2),$$
we can take a sequence of test functions $\vp_n$ such that, as $n\to\infty$,
$$\vp_n~\to~{\phi^\epsilon\over c(u)}\qquad \hbox{in}~~H^1(\R^2).$$
 Taking the limit, one obtains
\beq
\iint R\left[ \Big( {\phi^\epsilon\over c}\Big)_t - \phi^\epsilon_x\right]\, dxdt~=~
-\iint \frac{c'}{4c^2}(R^2-S^2)\phi^\epsilon \,dxdt.
\eeq
By (\ref{2.2}), this yields
\beq
\iint \left[{R\over c}(  \phi^\epsilon_t - c\phi^\epsilon_x) -{c'\over 2c^2}(R^2+RS)\phi^\ve
\right]\, dxdt~=~-\iint \frac{c'}{4c^2}(R^2-S^2)\phi^\epsilon \,dxdt\,,\eeq
and finally
\bel{R+S}
\iint {R\over c}(  \phi^\epsilon_t - c\phi^\epsilon_x)\,dxdt
~=~\iint\frac{c'}{4c^2}(R+S)^2\phi^\epsilon \,dxdt\,.
\eeq

By the way the test function $\phi^\epsilon$ has been defined at
(\ref{pen11})--(\ref{testf2})
the function
$\phi^\epsilon_t-c\phi^\epsilon_x$ is supported on a small neighborhood of
the boundary of $\Lambda(\Q)$.
More precisely, consider the four sets (Fig.~\ref{f:hyp203}, left)
\beq
\bega{rl}
\Gamma^\epsilon_{12}&\doteq ~\Big\{(t,x)\,;~~x_1^+(t)\leq x \leq x_1^+(t) +\epsilon \,,
~~x_1^-(t)-\epsilon\leq x \leq x_2^-(t) +\epsilon\Big\},\cr\cr
\Gamma^\epsilon_{13}&\doteq ~\Big\{(t,x)\,;~~x_1^-(t)-\epsilon\leq x \leq x_1^-(t)\,,
~~x_1^+(t)-\epsilon\leq x \leq x_2^+(t) +\epsilon\Big\},\cr\cr
\Gamma^\epsilon_{24}&\doteq ~\Big\{(t,x)\,;~~x_2^-(t)\leq x \leq x_2^-(t) +\epsilon \,,
~~x_1^-(t)-\epsilon\leq x \leq x_2^-(t) +\epsilon\Big\},\cr\cr
\Gamma^\epsilon_{34}&\doteq ~\Big\{(t,x)\,;~~x_2^+(t)-\epsilon \leq x \leq x_2^+(t)  \,,
~~x_1^-(t)-\epsilon\leq x \leq x_2^-(t) +\epsilon\Big\}.\cr\cr
\enda
\eeq
Notice that these sets  overlap near the points $P_i= (t_i, x_i)$, $i=1,2,3,4$.
However, each of these intersections is contained in a ball of radius $\O(\ve)$.
For example,
$$\Gamma_{12}^\epsilon\cap \Gamma_{13}^\epsilon~\subset~B(P_1, K\epsilon),
$$
for some constant $K$ and all $\epsilon>0$.
Observing that $R\in \L^2_{loc}(\R^2)$,  and $\|\phi^\epsilon_t + c(u) \phi^\epsilon_x\|_{\L^\infty}
= \O(\epsilon)$, we obtain the estimate
\beq\bega{l}
\ds\left|\iint_{\Gamma_{12}^\epsilon\cap\Gamma_{13}^\epsilon}\frac R{c(u)}(\phi_t^\epsilon-c\phi^\epsilon_x)\,dxdt\right|~
 \leq ~\ds\frac{ C_0}{ \epsilon} \iint_{B(P_1, K\epsilon)}|R|\, dxdt
\cr\cr
\qquad \leq ~\ds \frac{ C_0}{ \epsilon}
\int_{t_1-K\epsilon}^{t_1+K\epsilon}
\left(\int_{x_1-2K\epsilon}^{x_1+2K\epsilon}  R^2(t,x) dx\right)^{1/2}
(2K\epsilon)^{1/2} \, dt~\leq ~  \frac{ C_0}{ \epsilon}
E_0^{1/2}(2K \epsilon)^{3/2},
\enda
\eeq
for suitable constants $C_0, K$, and all $\epsilon >0$.
Repeating this argument for the other three intersections, we thus conclude
\beq
\bega{ll}
\lim\limits_{\epsilon\to 0}\iint {R\over c(u)}({\phi^\epsilon_t-c\phi^\epsilon_x})\,
dxdt~&\ds=\lim\limits_{\epsilon\to 0}\iint_{\Gamma^\epsilon_{12}\cup\Gamma^\epsilon_{13}\cup\Gamma^\epsilon_{24}\cup\Gamma^\epsilon_{34}}{R\over c(u)}({\phi^\epsilon_t-c\phi^\epsilon_x})\,dxdt\,\cr\cr
~&\ds=\lim\limits_{\epsilon\to 0}\left(\iint_{\Gamma^\epsilon_{12}}+\iint_{\Gamma^\epsilon_{13}}+\iint_{
\Gamma^\epsilon_{24}}+\iint_{\Gamma^\epsilon_{34}}\right)
{R\over c(u)}({\phi^\epsilon_t-c\phi^\epsilon_x})\,dxdt\,.
\enda
\eeq
Since $c(u)$ is  uniformly positive and bounded,
the same argument used in \eqref{h2} yields
\bel{iGe}
\lim_{\epsilon\to 0}~\iint_{\Gamma_{13}^\epsilon\cup\Gamma_{24}^\epsilon}
{R\over c(u)}({\phi^\epsilon_t-c\phi^\epsilon_x})\,dxdt ~=~0\,.\eeq

Concerning the integral over $\Gamma_{34}^\epsilon$, by (\ref{testf2}) we obtain
\beq
\bega{l}
\ds\lim_{\epsilon\to 0}\iint_{\Gamma_{34}^\epsilon}\frac{R(t,x)}{c(u(t,x))}
\cdot \frac{c(u(t,x_2^+(t)))+c(u(t,x))}{\epsilon}\, dxdt\cr\cr
\quad =~\ds \lim\limits_{\epsilon\to0}\iint_{\Gamma_{34}^\epsilon}  \frac{R(t,x)}{c(u(t,x))}
\frac{2c(u(t,x))}{\epsilon} dxdt+\lim_{\epsilon\to 0}\iint_{\Gamma_{34}^\epsilon}
\frac{R(t,x)}{c(u(t,x))}  \,\frac{c(u(t,x^+_2(t)))-c(u(t,x))}{\epsilon}\, dxdt\cr\cr
\quad =~\ds \lim_{\epsilon\to0}\frac2\epsilon\iint_{\Gamma_{34}^\epsilon}R(t,x)\,dxdt\,.
\enda
\eeq
Since $(R+S)^2\in \L^1_{loc}(\R^2)$, one has
\bel{lsource}
\lim_{\epsilon\to 0}\iint\frac{c'}{4c^2}(R+S)^2\phi^\epsilon \,dxdt
~=~\iint_{\Lambda(\Q)} \frac{c'}{4c^2}(R+S)^2 dxdt\,.\eeq
We thus conclude
\bel{lime2}
\lim_{\epsilon\to0}\frac1\epsilon\iint_{\Gamma_{34}^\epsilon}R\,dxdt-
\lim_{\epsilon\to0}\frac1\epsilon\iint_{\Gamma_{12}^\epsilon}R\,dxdt
~=~\ds\iint_{\Lambda(\Q)} \frac{c'}{8c^2}(R+S)^2 dxdt\,.
\eeq

Next, since $u\in H^1_{loc}$, we can write (Fig.~\ref{f:hyp203}, right)
\bel{lep3}\bega{rl}
u(P_4)-u(P_3)&\ds=~\lim_{\ve\to 0} {1\over\epsilon}
 \left( \int_{x_4-\epsilon}^{x_4}
 u(t_4, y)\, dy -
 \int_{x_3-\epsilon}^{x_3} u(t_3, y)\, dy\right)\cr\cr
 &\ds=~\lim_{\epsilon\to 0}
{1\over\epsilon}\dint_{\Gamma_{34}} \Big[u_t + c(u(t, x_1^+(t)))u_x
\Big]\, dxdt\cr\cr
&\ds=~
\lim_{\epsilon\to 0}{1\over\epsilon}
\dint_{\Gamma_{34}} \Big[u_t + c(u(t, x)u_x\Big]\, dxdt\cr\cr
&\ds=~
\lim_{\epsilon\to 0}{1\over\epsilon}
\dint_{\Gamma_{34}} R\, dxdt\,.
\enda
\eeq
This, and a similar estimate for $u(P_2)-u(P_1)$, yield
\bel{es5}
[u(P_4)-u(P_3)]-[u(P_2)-u(P_1)]~=~\iint_{\Lambda(\Q)}\frac{c'}{8c^2}(R+S)^2\, dxdt\,.
\eeq
Writing the right hand side of (\ref{es5}) as an integral
w.r.t.~the variables
$X,Y$,
by Remark~3
\beq
\bega{l}
[u(X_2, Y_2) - u(X_1, Y_2)] - [u(X_2, Y_1) - u(X_1, Y_1)]~=~\ds \iint_{\Lambda(\Q) }
\frac{c'}{8c^2}(R+S)^2\, dxdt\cr\cr
\qquad =~
\ds \iint_{\Q\cap \G } \frac{c'}{8c^2}(R+S)^2\cdot \frac{p}{(1+R^2)}\frac{q}{2c(1+S^2)}
\,
dXdY,\enda\eeq
where the last equality follows from Lemma 6(ii).
By Lemma  5, this shows that the weak derivative $u_{XY}$
exists and is given by
$$u_{XY}(X,Y)~=~\left\{\bega{cr}\ds -\frac{c'}{8c^2}(R+S)^2
\cdot\frac{p}{(1+R^2)}\frac{q}{2c(1+S^2)}
\qquad &\hbox{if}\qquad \det D\Lambda(X,Y) >0,\cr\cr
0\qquad &\hbox{if}\qquad \det D\Lambda(X,Y) =0.\enda
\right.$$
Recalling the definition (\ref{tvt}), we can write this weak derivative as
\bel{xpY}
u_{XY}~=~-\frac{c'}{16c^3}(\eta+\nu)pq+\frac{c'}{8c^3}(-\xi\zeta+\eta\nu)pq\,.
\eeq
By Lemma~4, it follows that, for a.e.~$X$, the map
$Y\mapsto  u_X(X,Y)~=~(\frac1{2c}\xi p)(X,Y)$ is absolutely continuous and its derivative is
given by (\ref{xpY}).
Recalling the equations for $p$ and $u$, since $p$ remains uniformly positive on bounded sets,
we conclude
\beq
\xi_Y~=~-\frac{c'}{8c^2}(\eta+\nu)q+\frac{c'}{4c^2}(\xi^2+\eta\nu)q\,.
\eeq
Similarly,
\beq
\zeta_X~=~-\frac{c'}{8c^2}(\eta+\nu)p+\frac{c'}{4c^2}(\zeta^2+\eta\nu)p\,.
\eeq
This completes the proof of Theorem~4. \endproof
\v

\section{Uniqueness of conservative solutions}
\setcounter{equation}{0}

We can now give a proof of Theorem~3, showing that conservative
solutions to the variational wave equation (\ref{1.1}) are unique.

Let initial data $u_0\in H^1(\R)$, $u_1\in\L^2(\R)$ be given.
These data uniquely determine a curve $\gamma$
in the $X$-$Y$ plane, parameterized by
$$X(x)~\doteq ~x+ \int_{-\infty}^x R^2(0,y)\, dy\,,\qquad\qquad
Y(x)~\doteq ~x+ \int_{-\infty}^x S^2(0,y)\, dy.$$
Along $\gamma$, the values of the variables
$(u,x,t,p,q,\nu,\eta,\xi,\zeta)$ are all determined by the data $u_0,u_1$.
Indeed, at the point $(X(x), Y(x))\in\gamma$ we have
$$\left\{\bega{rl}  t&=0  ,\cr\cr x&=~x\,,\enda\right.\qquad\qquad
\left\{\bega{rl}  u&=~u_0(x),\cr\cr p&=~q~=~1\,,\enda\right.
$$
$$\left\{\bega{rl}
\nu&\ds=~{1\over 1+R^2(0,x)}\,,\cr\cr   \eta&\ds=~{1\over 1+S^2(0,x)}\,,
\enda\right.
\qquad\qquad \left\{\bega{rl}\xi&\ds=~\, {R(0,x)\over 1+R^2(0,x)}\,,\cr\cr
 \zeta&=~\ds
\, {S(0,x)\over 1+S^2(0,x)}\,.\enda\right.$$
We recall that, by (\ref{2.1}),
$$R(0,x)~=~u_1(x) + c(u_0(x))\, u_{0,x}\,,\qquad\qquad
S(0,x)~=~u_1(x) - c(u_0(x))\, u_{0,x}\,.$$

Since the right hand sides of the equations in (\ref{sls}) are smooth,
given the above boundary data along $\gamma$, this semilinear
system has a unique solution in the $X$-$Y$ plane.  In particular, the functions
$(X,Y)~\mapsto~(x,t,u)(X,Y)$
are  uniquely determined, up to a set of zero measure in the $X$-$Y$ plane.
Since the map  $(x,t)\mapsto u(x,t)$ is continuous,
we conclude that $u$ is  uniquely determined, pointwise in the $x$-$t$ plane.
This completes the proof of  Theorem~3.
\endproof

{\bf Acknowledgment.} This research was partially supported by NSF, with grant DMS-1411786: Hyperbolic Conservation Laws and Applications, and the AMS Simons Travel Grants Program.


\end{document}